\documentclass[11pt]{iopart}

\usepackage{iopams}
\usepackage{latexsym}
\usepackage{graphicx}
\usepackage[export]{adjustbox}
\usepackage{color}
\usepackage{cite}
\usepackage{pgf,tikz}
\usepackage{hyperref}
\hypersetup{
     colorlinks = true,
     linkcolor = green!50!blue,
     anchorcolor = blue,
     citecolor = blue,
     filecolor = blue,
     urlcolor = blue
     }

\newcommand{\Gn}{\ensuremath{\boldsymbol \Gamma_{\!\rm{n}}}}
\newcommand{\Gpr}{\ensuremath{\boldsymbol \Gamma_{\!\rm{pr}}}}
\newcommand{\Gpo}{\ensuremath{\boldsymbol \Gamma_{\!\rm{pt}}}}
\newcommand{\dn}{\ensuremath{\mathbf d}}
\newcommand{\dnj}{\ensuremath{d_{j}}} 
\newcommand{\nnu}{\ensuremath{\boldsymbol \nu}}
\newcommand{\nnumap}{\ensuremath{\nnu_{\!\rm{\scriptscriptstyle MAP}}}}

\newcommand{\proofend}{\hfill $\Box$}

\begin{document}

\title[Bayesian approach to inverse scattering with topological priors]
{Bayesian approach to inverse scattering with topological priors}

\author{Ana Carpio$^{1,2}$, Sergei Iakunin$^{1,3}$, Georg Stadler$^2$}

\address{$^1$ Universidad Complutense de Madrid, 28040 Madrid, Spain}

\address{$^2$ Courant Institute, New York University, NY 10012, USA}

\address{$^3$ Basque Center for Applied Mathematics-BCAM, 48009 Bilbao, Spain}

\eads{\mailto{ana\_carpio@mat.ucm.es},
\mailto{siakunin@ucm.es}, \mailto{stadler@cims.nyu.edu}}

\begin{abstract} 
We propose a Bayesian inference framework to estimate uncertainties in inverse
scattering problems. Given the observed data, the forward model and
their uncertainties, we find the posterior distribution over a finite
parameter field representing the objects. To construct the prior
distribution we use a topological sensitivity analysis.  We
demonstrate the approach on the Bayesian solution of 2D inverse
problems in light and acoustic holography with synthetic
data. Statistical information on objects such as their center
location, diameter size, orientation, as well as material properties,
are extracted by sampling the posterior distribution.  Assuming the
number of objects known, comparison of the results obtained by Markov
Chain Monte Carlo sampling and by sampling a Gaussian distribution
found by linearization about the maximum a posteriori estimate show
reasonable agreement. The latter procedure has low computational cost,
which makes it an interesting tool for uncertainty studies in
3D. However, MCMC sampling provides a more complete picture of the
posterior distribution and yields multi-modal
posterior distributions for problems with larger measurement noise. 
When the number of objects is unknown, we devise a stochastic model 
selection framework.  
\end{abstract}




\section{Introduction}
\label{sec:basic}

Inverse scattering techniques are a common tool to detect objects
in areas such as medicine, geophysics, or public security. The basic
structure of the underlying mathematical problem is as follows.
An incident wave field illuminates a set of objects integrated in an
ambient medium. The resulting wave field is measured at a set of 
detectors. Given the measured data, the goal is to reconstruct the
unknown objects and their material properties. In practice, the
process is affected by different sources of errors and uncertainty,
such as external noise in the recorded data and errors in the 
measurement systems and governing mathematical
equations. 
Some deterministic approaches are able to provide reasonable
reconstructions of objects  under specific conditions.
However, these reconstructions depend on the choice of 
tuning parameters, such as thresholds in direct methods and 
regularization or stopping criteria in iterative procedures.
Moreover, deterministic approaches do not provide
information on the confidence we are allowed to have in the results
and do not shine light on correlations between the inferred parameters. 
This gap is addressed by a Bayesian probabilistic formulation which 
provides a more complete picture of the reconstructed parameters 
and their uncertainties.

Details of the imaging process depend on the type of waves employed. To 
fix ideas, we focus here on situations where the physical process is
modelled by a wave equation.  This is the case, for instance, in inverse acoustic 
scattering \cite{buishape,topologicalenergy,davidacoustic,georgmarmousi}, 
and inverse  electromagnetic  scattering when using polarized radiation
\cite{belkebir,jcp19,dorn,litman}. When the incident wave is 
time harmonic,  $U_{\rm inc}(\mathbf x,t) = e^{-\imath \omega t} 
u_{\rm inc}(\mathbf x)$, the total wave  field is time harmonic too, i.e., 
$U(\mathbf x,t) = e^{-\imath \omega t} u(\mathbf x)$. 
Its amplitude $u(\mathbf x)$ obeys a Helmholtz transmission 
problem, which, in two dimensions, is
\begin{equation}\label{forward}
\hskip -1cm\left\{\begin{array}{ll}
\Delta u + \kappa_{\rm e}^2 u = 0   &
\mbox{\rm in $\Omega_{\rm e}$},  \\[1ex]
\Delta u + \kappa_{\rm i}^2 u =0  &
\mbox{\rm in $\Omega_{\rm i}$},  \\[1ex]
u^{-} - u^{+} = 0,\quad   
\beta {\partial u^{-} \over \partial \mathbf n} - 
{\partial u^{+}  \over \partial \mathbf n}= 0
&\mbox{\rm on $\partial \Omega_{\rm i}$},  \\[1ex]
{\rm lim}_{|\mathbf x| \rightarrow 0} |\mathbf x|^{1/2} 
\left({\partial \over \partial |\mathbf x|}(u-u_{\rm inc}) 
- \imath k_{\rm e} (u-u_{\rm inc}) \right) = 0,
\end{array}
\right.
\end{equation}
where $\Omega_{\rm i}$ is an inclusion, $\Omega_{\rm e}= \mathbb R^2
\setminus {\overline{\Omega}_{\rm i}}$, and $\mathbf n$ is the unit
outer normal vector for $\Omega_{\rm i}$.
The symbols $^-$ and $^+$ denote values from inside
and outside $\Omega_{\rm i}$, respectively. The Sommerfeld radiation
condition on the propagation of the scattered field $u_{\rm
  sc}=u-u_{\rm inc}$ at infinity implies that only outgoing waves are
allowed.  The parameters $\beta, \kappa_{\rm
  i},\kappa_{\rm e}\ge 0$ depend on the frequency $\omega$
and the material properties.  Moreover, $\kappa_{\rm e} \sim k_{\rm
  e}$ is assumed to be constant outside of a ball  containing the
objects and the detectors $ \mathbf x_j,$ $j=1,\ldots,N$. 

\begin{figure} \centering
\includegraphics[width=11cm]{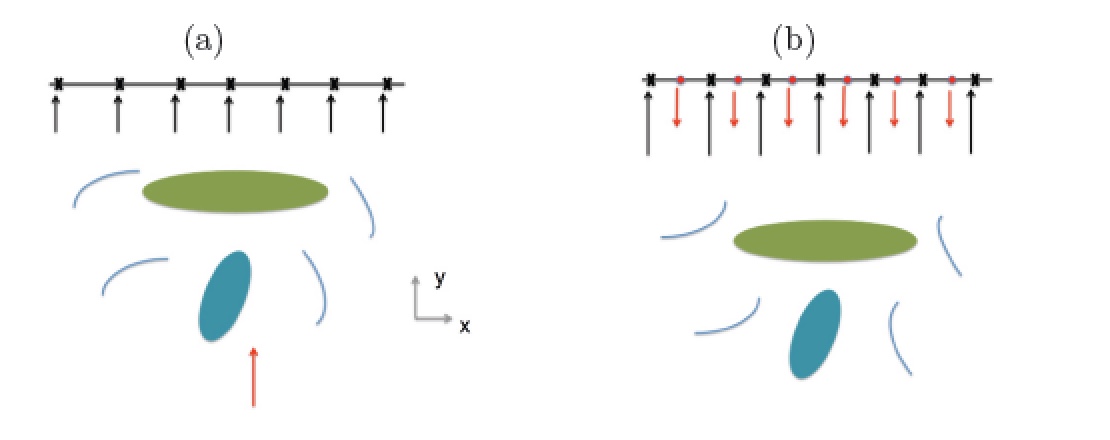}
\caption{Schematic arrangements of incident waves (red beams) and 
detector grids (black crosses) in different imaging set-ups. 
(a) In holography and microscopy, one incident planar wave interacts 
with the objects and is recorded at a fixed grid of detectors behind the 
objects \cite{davidacoustic,davidlight}. 
(b) In acoustic imaging of materials, waves emitted from a grid of 
sources interact with the medium and the reflected waves 
are recorded at a grid of  receivers in the same region
\cite{greengard,feijoooberai,guzinah}.
Blue lines represent waves interacting with the objects
and black arrows indicate the reflections arriving at
the detector grid.
Microwave imaging  uses set-ups similar to (a) for multiple incident wave 
frequencies, while the relative positions of objects, emitters and receptors 
are rotated to increase the number of independent observations \cite{belkebir,litman}.} 
\label{fig1}
\end{figure}

To model the measurement process, we generate data $\dn$ by solving
the forward problem (\ref{forward}) with inclusions we consider as the
``truth'' and evaluate the resulting wave field at detectors located at
$\mathbf x_j$, $j=1,\ldots,N$. We then add independent additive white noise of
a magnitude specified in each problem  to obtain ``synthetic'' measurement.
Depending on the application, the measured data are complex-valued
amplitude fields $u(\mathbf x_j)$ (in microwave imaging
\cite{belkebir,litman} or acoustic holography \cite{davidacoustic},
for instance) or real-valued intensities $|u(\mathbf x_j)|^2$ (in light
microscopy \cite{davidlight, tombayesian}).  Given the measured data
and the ambient properties $\kappa_{\rm e}$, our goal is to find the
inclusions $\Omega_{\rm i} = \cup_{\ell=1}^L\Omega_{\rm i}^\ell$.
Different experimental imaging set-ups correspond to different
arrangements of emitted incident waves and detector distributions, see
Figure \ref{fig1}.  Here, we focus on the configuration displayed
in Fig.~\ref{fig1}(a) for our numerical tests, though the methods
extend to other arrangements. We will consider both types of data,
complex valued fields and real intensities, which is the case in
acoustic and light holography
\cite{tombayesian,davidacoustic,davidlight}, respectively. Bayesian
methods are particularly interesting for acoustic imaging, since the
magnitude of noise in the recorded data is usually larger.

If we assume that each object can be approximated by 
star-shaped parameterizations, that is, the boundary of each $\Omega_{\rm i}^\ell$ 
is defined by `rays' emerging from a center, the inverse scattering 
problem amounts to finding a set of parameters $\nnu$: 
the number of components,  their centers and the coefficients of 
the trigonometric expansions defining their boundaries, as well as 
their material  constants  $\kappa_{\rm i}^\ell$.
The coefficients $\kappa_{\rm i}^\ell$ may be spatially dependent, which increases 
the complexity of the problem. Here, we  take them to be 
constant and  known. We represent star-shaped 
objects using expansions in trigonometric polynomials for the radius 
\cite{ip08,jcp19,hohage2d} instead of general functions varying with the 
angle  \cite{buishape,dunlopthesis} to reduce the parameter dimension.

There is a broad literature on deterministic inverse scattering 
approaches, which basically fall in two categories. Non-iterative methods 
such as linear sampling, factorization, multiple signal classification, 
orthogonality sampling, direct sampling and topological derivative analysis 
\cite{cakonilsm, kirschmusic, dsm, feijoo} 
provide direct approximations to the objects from the data. 
Iterative techniques may refine this information at a higher 
computational cost, and vary widely with the specific application.
In most set-ups, one can resort to regularized nonlinear least-squares 
optimization formulations, where the governing equation is treated as 
an equality constraint.
The idea is to seek shapes that fit the measurements within the noise
level, and to use regularization or early termination of iterative
algorithms to prevent overfitting. 
In principle, the existence of several local minima is a
possibility---if local minima occur depends on the distribution of
detectors, the incident wave directions and the local wave speed. To
obtain the best estimate for the true configurations, one may follow
different strategies: produce sharp initial guesses of the objects or
select wide ranges of detectors and incident waves or frequencies
\cite{ip08,conca,osher,park,guzinam,laurain}, when possible,
and include additional regularizing terms in the cost functional, see
\cite{jcp19,hohage2d} for instance.

Such optimization-based methods are naturally related to Bayesian
approaches, which consider all variables in the inverse problem as
random variables. Assuming the variables are
a finite-dimensional vector $\nnu$, the densities of the random variables are
related using Bayes' formula \cite{KaipioSomersalo06,Tarantola05}:
\begin{eqnarray}
p_{\rm pt}(\nnu):= p(\nnu | \dn) = { p( \dn | 
\nnu) \over p(\dn)} p_{\rm pr}(\nnu).
\label{bayes}
\end{eqnarray}
Here, $p_{\rm pr}(\cdot)$ is the prior density of the variables, which
incorporates prior (expert) knowledge available about the variables;
$p(\dn|\nnu)$ is the likelihood or conditional probability
of the observations $\dn$ given the variables $\nnu$; and $p_{\rm
  pt}(\dn|\nnu)$ is the posterior density of the
parameters given the data, which is the solution of the Bayesian
inverse problem. The density $p(\dn)$ is a scaling that does not
depend on the parameters. Full characterization of the posterior
density is an extremely challenging probability problem for
moderate and high-dimensional parameters $\nnu$, and often
one has to rely on approximations of the posterior distribution. One
approximation is finding the maximum a posteriori (MAP) point, i.e.,
the parameters that maximize the posterior density. This amounts to an
optimization problem that, making assumptions on the prior and
likelihood and taking the negative logarithm of the densities, becomes
a nonlinear least-squares problem of the form used in
deterministic inversion;
see Section \ref{sec:bayesian}. In this optimization problem, the
regularization is implied by the prior density and the likelihood
corresponds to the data misfit term. Linearization about the
MAP point, also called Laplace approximation, results in an
approximation of the posterior density by a multivariate Gaussian, 
allowing for computationally efficient manipulations even for high- or
infinite-dimensional parameters \cite{georglinearized,georgmarmousi}.
Alternatively, to avoid a Gaussian approximation of the posterior, we
can sample the posterior distribution using, for instance, Markov
Chain Monte Carlo (MCMC) techniques
\cite{dunlopthesis,goodmanweare,georgmcmc} at a higher computational
cost.

The results obtained with Bayesian methods depend on the choice of
the prior \cite{tombayesian,dunlopthesis}. Prior distributions can
involve additional parameters, so-called `hyperparameters'
\cite{hyperparameters}.  There are different approaches to handling
hyperparameters. A straighforward strategy consists in fixing
subjective guesses for them. The results depend on how good such
subjective guesses are.  A second possibility is to select empirical
guesses, i.e., one solves the problem using different hyperparameters
and then chooses the most appropriate one
\cite{empiricalbayes}. Finally, one might introduce additional
probability densities for the hyperparameters, and work in an
hierarchical Bayes framework \cite{hierarchicalbayes}.  In this work,
we construct subjective guesses for the prior distributions from a
topological sensitivity analysis of the cost functionals.
This allows us to handle imaging set-ups with rather limited data, i.e.,
in which detectors are located in a narrow area, only one incident wave 
is used, multiple penetrable objects may be present,
and the recorded data may not be the complex amplitude field, 
but just its real modulus.
We will show that both the Bayesian linearized approach and MCMC 
sampling provide reasonable descriptions of statistical properties of the
objects for various noise levels in the observations, when the number
of object components is known. 
Previous work in 3D light holographic settings
also applied MCMC sampling to infer the location, size and refractive index of single 
spherical particles \cite{tombayesian}. In
\cite{palafox1,palafox2,palafox3}, scattering from single 2D
sound-soft objects is considered. Here, the objects are
placed at known locations and inference is based on far-field,
small-noise, complex data. Similar to this work, MCMC sampling is used
to infer Fourier modes in a starshaped
parametrization, amongst other parametrizations.
Bayesian approaches relying on more
complex parameterizations varying with the angle are presented in
\cite{buishape,dunlopthesis} for 2D acoustic and tomography set-ups.

When the number of objects components is not known, this number
becomes an additional unknown. Compared to the other parameters, it is
of a rather different nature: it is discrete, and it controls the
presence of other parameter blocks. In this case, we reformulate our
problem within the frameworks of hierarchical Bayesian modeling and
model selection \cite{modelselection}.  The number of objects can then
be selected by empirical arguments or by MCMC methods with variable
selection \cite{vbleMCMC}, see also 
\cite{dunloplevel} for a hierarchical level set approach.
\cite{palafox1} applies model selection to estimate the number
of modes in the radii of single starshaped objects instead.

The next sections are organized as follows. In
Section~\ref{sec:bayesian}, we devise a general Bayesian approach for
object reconstruction, which we test in the physical set-ups described 
in Section~\ref{sec:physical}. Section~\ref{sec:priors} uses
topological sensitivity analysis to estimate 
the number of objects and define priors. Assuming the number
of objects and their material properties known, we compute MAP points
for their shape and then sample the Gaussian posterior approximation
obtained by linearization at the MAP point in Section \ref{sec:bl}. The 
results for different geometries are compared to MCMC sampling of the 
full posterior distribution in Section~\ref{sec:bmcmc}. 
Section~\ref{sec:material} discusses cases where we allow for variations 
in the material constants of the inclusions. In Section~\ref{sec:model},
the number of objects is considered as unknown, and thus we
implement a Bayesian strategy for model selection.
Finally, Section~\ref{sec:conclusions} presents our conclusions.  
For completeness, a final Appendix
discusses formulas for different derivatives relevant in
Sections~\ref{sec:bl} and \ref{sec:material}. 

\section{Bayesian approach for a given number of objects}
\label{sec:bayesian} 

In this section we develop a framework for object detection. We assume
that we can approximate the objects using star-shaped parameterizations,
but more general parameterizations could be used instead. The
parameters for an inclusion consisting of $L\ge 1$ object components
are collected in a vector $\nnu:=(\nnu^1,\ldots,\nnu^L)\in \mathbb R^{L(2M+3)}$, 
where
\begin{equation} 
\nnu^\ell = (c_{x}^\ell, c_{y}^\ell,a_{0}^\ell,a_{1}^\ell,\ldots,
a_{M}^\ell,b_{1}^\ell,\ldots,b_{M}^\ell),  \quad \ell=1,\ldots,L. \label{parameters}
\end{equation}  
Here, $(c_x^\ell,c_y^\ell)$ are the centers and $r^\ell(t)$ the radii of the 
object components, associated to the parameterization 
\begin{eqnarray}
\label{trigonometric1}
\mathbf q^\ell(t)=
(c_x^\ell,c_y^\ell)+r^\ell(t)(\cos(2 \pi t),\sin(2 \pi t)),\quad t \in[0, 1],\\
\label{trigonometric2}
r^\ell(t)= a_0^\ell + 2 \sum_{m=1}^M a_m^\ell \cos(2 \pi m t) +
 2 \sum_{m=1}^M b_m^\ell \sin(2 \pi m t),
\end{eqnarray}
for $\ell=1,\ldots,L$. 
Note that similar parametrizations are available in three
dimensions, where they could rely for instance on spherical harmonics instead of Fourier
series \cite{jcp19}.
The number of modes $M$, which we fix, controls
the complexity of the boundary.  The Bayesian approach to inverse
problems requires to define a prior distribution for the parameters
$\nnu$. We choose $p(\nnu)$ as a multivariate Gaussian
\begin{eqnarray} \begin{array}{l}
p(\nnu)  
={1\over (2 \pi)^{n/2}}  {1\over \sqrt{|\Gpr |}}
\exp(-{1\over 2}
(\nnu - \nnu_0)^t \Gpr^{-1} 
(\nnu - \nnu_0) ) 
\end{array} \label{prior}
\end{eqnarray}
with covariance matrix $\Gpr$ and $n:= L(2M+3)$. Notice that the radii
(\ref{trigonometric2}) belong to the space of trigonometric
polynomials $T_{2M+1},$ and are expanded in an orthonormal basis
$\phi_1$,\ldots,$\phi_{2M+1}$, so that the mass matrix associated to
this basis is the identity. Otherwise, the mass matrix would enter
(\ref{prior}), see \cite{georglinearized}. Modeling of the prior
distribution, i.e., the choices for the covariance matrix $\Gpr$ and the 
mean $\nnu_0$, are discussed in Section \ref{sec:priors}.

We illuminate the $L$ objects with an incident plane wave of amplitude
$u_{\rm inc}(\mathbf x)$, generating data $\dn$ at detectors $\mathbf
x_j$, $j=1,\ldots, N$. We denote by $\mathbf f:\mathbb
R^{L(2M+3)}\rightarrow \mathbb R^N$ the parameter-to-observable map,
i.e., the mathematical description of this process. To be more
precise, for parameters $\nnu$, we denote by $u_{\Omega_{\nnu}}$ the
solution of the wave equation (\ref{forward}) with object $\Omega_{\rm
  i}= \Omega_{\nnu}$ defined by
(\ref{parameters})-(\ref{trigonometric2}). Then $\mathbf f
(\nnu)\:=\left( f(u_{\Omega_{\nnu}}(\mathbf x_j)) \right)_{j=1}^N$,
where $f$ is the measurement operator (which may be real $f(u)= |u|^2$
or complex valued $f(u)=u$).

We assume additive Gaussian measurement noise, i.e., the observations
and parameters are related by
\begin{equation}\label{eq:IP}
  \dn = \mathbf f(\nnu) + \boldsymbol\varepsilon.
\end{equation}
Here, the measurement noise $\boldsymbol\varepsilon$ is distributed as a
multivariate Gaussian ${\cal N}(0,\Gn)$ with mean
zero and covariance matrix $\Gn$. We consider the noise level for each
sensor to be equal and uncorrelated, so that $\Gn$ is a real
diagonal matrix.  For complex valued data, the additive noise is
represented by a standard complex Gaussian variable whose real and
imaginary parts are both real Gaussians of the form ${\cal
  N}(0,\Gn/2)$ \cite{complex}.

Due to these assumptions, on the measurement noise, the
conditional probability density $p( \dn | \nnu)$ takes the form
\begin{eqnarray}
p( \dn | \nnu) = {1 \over (2\pi)^{N/2} 
\sqrt{|\Gn|}} \exp \Big(- {1 \over 2} \|
\mathbf f( \nnu) - \dn  \|^2_{\Gn^{-1}} \Big),
\label{likelihood}
\end{eqnarray}
where $\| \mathbf v \|_{\Gn^{-1}}^2 =  
\mathbf {\overline v}^t \Gn^{-1} \mathbf v$.
Combining (\ref{bayes}), (\ref{likelihood}), (\ref{prior}) and neglecting
normalization constants, the posterior
density becomes, up to multiplicative constants,
\begin{eqnarray}
p_{\rm pt}(\nnu) \propto \exp \left( -{1\over 2} \|
\mathbf f( \nnu) - \dn  \|^2_{\Gn^{-1}} -
 {1\over 2} \|  \nnu - \nnu_0 \|_{\Gpr^{-1}}^2
 \right). \label{posterior}
\end{eqnarray}

Taking logarithms, the problem of maximizing the posterior probability
of the parameter set $\nnu$ given the data $\dn$ is identical to
minimizing the regularized cost objective \cite{optimizationbayes}:
\begin{eqnarray}
J(\nnu) := {1\over 2}\|\mathbf f( \nnu) 
- \dn \|^2_{\Gn^{-1}} 
+  {1\over 2} \|\nnu - \nnu_0\|^2_{\Gpr^{-1}},
\label{regcost}
\end{eqnarray}
where we neglect additional terms involving only the covariances.
The first part of functional (\ref{regcost}) is related to the standard
cost used in deterministic inverse problems,
whereas the second part originating from the prior takes the role of the
regularization, which prevents ill-posedness and overfitting of the
observation data.

In our experiments, we choose a diagonal covariance matrix $\Gn= {\rm
  diag}(\sigma_1^2,\ldots,\sigma_N^2)$, and set all the variances
equal to a constant $\sigma_{\rm noise}^2.$ Thus, $\sqrt{|\Gn|}=
\Pi_{j=1}^N \sigma_j = \sigma_{\rm noise}^N$. The mean $\nnu_0$ in the
prior multivariate Gaussian density and the elements of the covariance
matrix $\Gpr$ are considered hyperparameters to be selected as
discussed in the next section.  This selection reflects our
uncertainty in the prior information available, and affects the
resulting parameter inference. We introduce a strategy to generate guesses 
for all the hyperparameters based on the topological fields of the
objective functional (\ref{cost}). For that purpose, the available
data are split in two parts: a fraction $\dn^{(1)}$ is used to
generate the prior, whereas the remaining values $\dn^{(2)}=:\dn$ are
used in (\ref{posterior}) to define the posterior distribution.
We use intersperse grids as detailed in Section \ref{sec:blresults}.
We generate synthetic data $\dn^{(1)}$ and $\dn$ for our tests solving
(\ref{forward}) by BEM methods using `true' object inclusions, and
then add noise (as detailed in our numerical tests) to these
observations.

\section{Physical set-ups}
\label{sec:physical}

We will study the behavior of the methods presented in this paper in
light and acoustic holography set-ups, adjusted to Figure
\ref{fig1}(a). In this section we briefly summarize the physics
background and the parameter choices we make for the remainder of this
paper.

Let us fix a reference length scale $L$. Typical object sizes may range 
from $L/10$ to $2L$, for instance, while the distance to the detectors 
is about $5L$. 
In an acoustic setting, the evolution of the total wave field $U$
is governed by the wave equation $\rho(\mathbf x) U_{tt}(\mathbf x,t) -
{\rm div} (\alpha (\mathbf x) \nabla U(\mathbf x,t)) = 0 $,
where $\rho$ is the density and $\alpha$ represents a relevant  `modulus'.
In liquids,  $\alpha$ represents the bulk modulus and
in solids, $\alpha$ is Young's modulus. In gases, $\alpha = \gamma p$,
where $p$ is the pressure and $\gamma$ is related to the specific heat.
The speed of sound in the medium is then $c =  (\alpha/\rho)^{1/2}$. 
When the applied field is time-harmonic, then $U(\mathbf x, t) = 
e^{-\imath \omega t} u(\mathbf x)$, where the amplitude $u(\mathbf x)$ 
obeys  ${\rm div} (\alpha (\mathbf x) \nabla u(\mathbf x)) 
+\rho(\mathbf x) \omega^2 u(\mathbf x) = 0$.

We nondimensionalize the problem according to $x = x' \, L$, $y = y' \, L$, 
$\Omega = \Omega' L$, $u_{\rm inc} = u_0 u_{\rm inc}'$, $u= u_0 u'$ and then 
drop the $'$ for ease of notation. Here, $L$ is the reference length and 
$u_0$ the modulus of the amplitude of the incident wave. The
total amplitude field is then governed by (\ref{forward}), where the
parameters are the dimensionless wavenumbers,
\begin{eqnarray}
\kappa_{\rm e}=\frac{2\pi\nu }{c_{\rm e}} L, \quad
\kappa_{\rm i}=\frac{2\pi\nu }{c_{\rm i}} L, \label{wavenumbersound}
\end{eqnarray}
as well as the ratio $\beta = {\alpha_{\rm i} \over \alpha_{\rm e}}$. We denote
by $\nu={\omega \over 2 \pi}$ the frequency of the emitted sound wave and 
by $c_{\rm e}$ and $c_{\rm i}$ 
the sound speed in the ambient medium and inside the objects, 
respectively. The incident amplitude for the plane wave is
$u_{\rm inc}= e^{\imath \kappa_{\rm e} y}$, where $y$ points in the direction 
of the detectors. 

Frequencies for sound lie in the range $20$Hz--$20$kHz. The speed of 
sound in air  is $343$m/s.  If for instance $L \sim 6$cm, 
 $\kappa_{\rm e} \sim 12$-$20$ for $11$--$18$kHz. Depending on the
material the object is made of, it can be sound-soft (it absorbs sound),
sound-hard (it reflects sound), or penetrable. We assume
penetrable objects. For the other types, the governing equations
and formulas should be adjusted following \cite{ln08}.

In a light holography set-up, the general framework is similar \cite{sims18,jcp19}, 
assuming we use polarized light in the presence of few well separated objects. 
Then, the equations governing the amplitude field can be approximated by 
(\ref{forward}) and  the wavenumbers are again given by (\ref{wavenumbersound}), 
$c$ and $\nu$ representing light wavespeed and frequencies.   In this
case, $\beta ={\mu_{\rm e} \over \mu_{\rm i}}$, $\mu_{\rm e}$ and $\mu_{\rm i}$
are material permeabilities.
For the biological applications we target, $\beta \sim 1$.
Visible light wavelengths lie the range of $400-700$nm and result in 
wavenumbers $\kappa_{\rm e}\sim 12$-$20$, setting $L \sim 1\mu$m, for instance. 
While in classical light holography, the measured data are real
intensities $|u|^2$ \cite{sims18}, it has become possible to record
complex-valued data $u$ in acoustic holography set-ups
\cite{davidacoustic}. The latter are affected by larger noise
magnitudes.

\section{Topological selection of priors}
\label{sec:priors}

Topological derivative methods generate first guesses of objects without 
a priori information, other than the measured noisy data,  the ambient 
medium properties and the incident wave.  In deterministic frameworks, 
such guesses are then improved by level set, shape derivative, topological  
derivative or Gauss-Newton iterations \cite{ip08, jcp19, conca, osher}.
Following a Bayesian approach,  we propose the following procedure to generate 
prior densities from topological sensitivity studies 
of the underlying unregularized cost functional (\ref{cost}).
Such prior knowledge can significantly influence Bayesian inference
results, as, e.g., observed in \cite{tombayesian,dunlopthesis}.

\subsection{Topological derivative of the cost functional}
\label{sec:topologicalfields}

Given ${\cal R} \subset \mathbb R^2$, the topological derivative of  
\begin{eqnarray}
J_c(\mathbb R^2 \setminus \overline{\Omega_{\nnu}}) := {1\over 2}
\sum_{j=1}^N |\dnj - f(u_{\Omega_{\nnu}}(\mathbf x_j))|^2,
\label{cost}
\end{eqnarray}
is  a scalar field $D_{\rm T}(\mathbf{x}, {\cal R})$ satisfying \cite{sokolowski}:
\begin{eqnarray*} \label{expansion}
J_c({\cal R} \setminus \overline{B(\mathbf x, \varepsilon)}) = J_c( {\cal R} )
+ D_{\rm T}(\mathbf{x}, {\cal R}) \, {\rm meas}(B(\mathbf x, \varepsilon)) + o(\varepsilon^2)
\end{eqnarray*}
for any $\mathbf x \in {\cal R}$ and any small radius $\varepsilon >
0$. Here, ${\rm meas}(B(\mathbf x, \varepsilon))$ is
the volume of the ball centered at  $\mathbf x$ with radius
$\varepsilon.$
When $D_{\rm T}(\mathbf{x}, {\cal R}) <0$, the cost functional
decreases by removing small balls centered at $\mathbf x$. This suggests
that placing objects $\Omega$ in regions where the topological
derivative is negative, the cost functional should
decrease. In this way, we obtain guesses for object component locations.
When ${\cal R}= \mathbb R^2$, the topological
derivative admits an explicit expression in terms of auxiliary forward
and adjoint fields \cite{ip08,guzinah}.  For $\mathbf x \in \mathbb
R^2$
\begin{eqnarray} \label{TDempty}  \hskip -1cm
D_{\rm T}(\mathbf{x},\mathbb{R}^2) = \mbox{\rm Re}\left[
{2(1-\beta) \over 1 + \beta} \nabla u_{\rm inc}(\mathbf{x}) \nabla \overline{p}_{\rm inc}(\mathbf{x}) 
+ (\beta \kappa_{\rm i}^2 \!-\! \kappa_{\rm e}^2)\, u_{\rm inc}(\mathbf{x}) \overline{p}_{\rm inc}(\mathbf{x})
\right] , 
\end{eqnarray}
where $\kappa_{\rm e}$ is the wavenumber for the outer medium, $\kappa_{\rm i}$ the wavenumber for the inclusion whose shape we are seeking and $u_{\rm inc}$ is the incident 
wave. When $\kappa_{\rm e}$ is constant, the conjugate adjoint field is
\begin{eqnarray}
\overline p_{\rm inc}(\mathbf x) = {\imath \over 4} \sum_{j=1}^N
H_0^{(1)}(\kappa_{\rm e} |\mathbf x - \mathbf x_j|) \chi(\mathbf x_j),
\label{pinc}
\end{eqnarray}
with $H_0^{(1)}$ the Hankel function and $\chi(\mathbf x_j)= \overline{(\dnj -
f(u_{\rm inc}(\mathbf x_j) )} f'(u_{\rm inc}(\mathbf x_j)).$ Standard choices for the 
measurement operator are $f(u)=u$ with $f'(u)=1$ and $f(u) = |u|^2$ with 
$f'(u)= 2 \overline{u}.$ 
When $\kappa_{\rm e}$ varies spatially inside a bounded region, 
 $\overline{p}_{\rm inc}$ is a solution of:
\begin{eqnarray*} \label{adjointempty}
\begin{array}{l}
 \Delta \overline{p} + {\kappa_{\rm e}^2} \overline{p} =   \sum_{j=1}^N \chi(\mathbf x_j)
 \delta_{\mathbf x_j} \; \mbox{ in } \mathbb{R}^2, \quad
\lim_{|\mathbf x|\to \infty}|\mathbf x|^{1/2} \left( {\partial \overline{p} \over
\partial |\mathbf x| } - \imath k_{\rm e} \overline{p} \right)=0,
\end{array}
\end{eqnarray*}
$k_{\rm e}$ being its constant value at infinity.

\subsection{Prior selection}
\label{sec:prior}

\begin{figure}[h!]
\centering
\includegraphics[width=14cm]{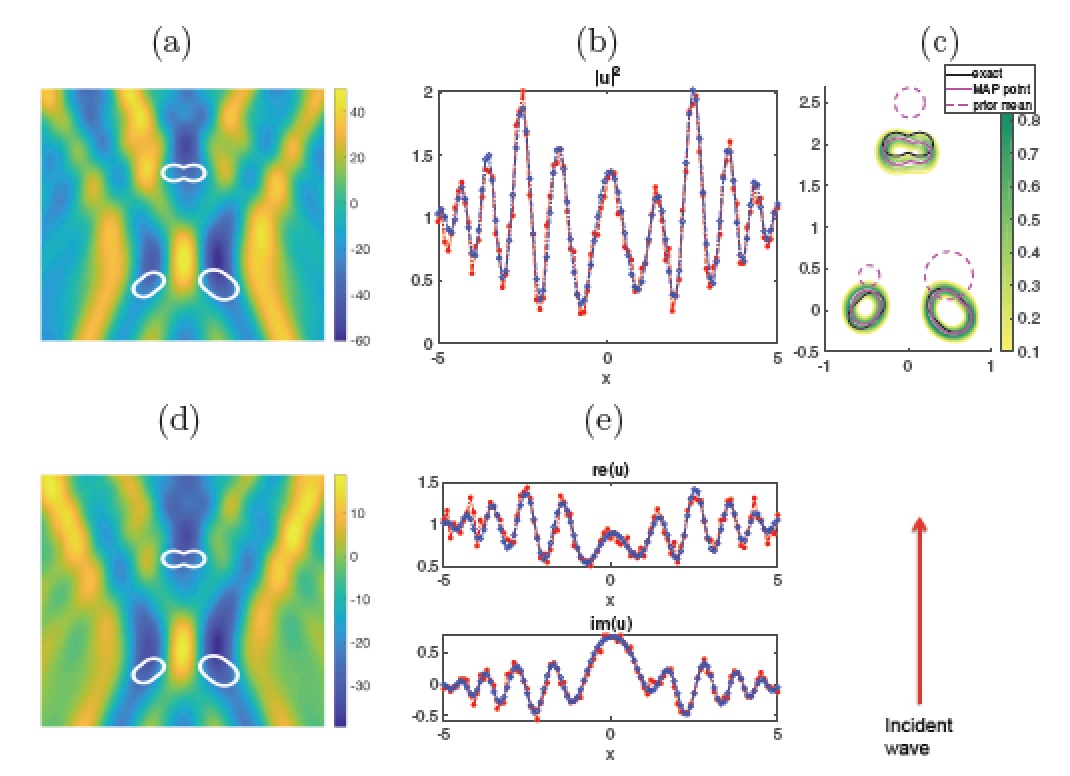}
\caption{Topological derivative (\ref{TDempty}) for a configuration
  with three objects whose contours are superimposed: (a) with data
  $|u(\mathbf x_j)|^2$, $j=1,\ldots,101$, (d) with data $u(\mathbf
  x_j)$, $j=1,\ldots,101$, corrupted by noise of magnitude $5\%$ and
  $10\%$, respectively.  The noisy data, depicted in panels (b) and
  (e), are recorded in the set-ups represented in Figure
  \ref{fig1}(a). Red asterisks represent the noisy data, whereas blue
  circles are the `true' synthetic values generated solving the
  forward problem by BEM. Panel (c) represents the information on the
  objects inferred for (a)-(b) by the Bayesian linearized approach
  developed in Section \ref{sec:bl} with $M=5$.  Parameter values:
  $\kappa_{\rm i}=15.12$, $\kappa_{\rm e}=12.56$, $\beta=1$. }
\label{fig2}
\end{figure}

Figure \ref{fig2} displays the topological  derivative field (\ref{TDempty})  
for three objects when $\kappa_{\rm i}=15.12$ and $\kappa_{\rm e}=12.56$.  Regions where large negative values are attained provide guesses of the object locations, which we use as prior knowledge as follows. We choose a constant $C_0 \in (0,1)$ and define the set
$\Omega_0 = \{ \mathbf x \in R_{\rm obs} | D_{\rm T}(\mathbf x, \mathbb R^2) < (1-C_0) 
{\rm min}_{\mathbf y \in R_{\rm obs}} D_{\rm T}(\mathbf y, \mathbb R^2) \}
$ in the region where we look for objects, the so-called observation region $R_{\rm obs}.$ Then, we fit circles to each connected component to construct our  initial guess $\Omega_{\nnu_0}.$

For single objects, we  test $50$ values between $0.01$ and
$0.3$ for $C_0$, choosing the one which yields the
smallest value for the cost functional. Here, depending on the value
selected for $C_0$, we will capture one, two or three dominant
negative regions. Since we assume we know the number of objects, we
set $C_0$ to capture $L=3$ objects.  Star-shaped objects are a
deformation of circles.  We fit circles to each component
$\Omega_{\nnu_0}^\ell$ as follows. For $\ell=1,\ldots, L$,
$(c_{x,0}^\ell,c_{y,0}^\ell)$ is the center of mass of the component
and $a_{0,0}^\ell$ is the minimum distance from the center to the
boundary (or the average of minimum and maximum distances).  We set
$a_{m,0}^\ell=0$ and $b_{m,0}^\ell=0$, $m=1,\ldots,M$ to avoid
inserting unnecessary bias in the hyperparameters.  Finally, we check
that $J_c(\mathbb R^2 \setminus{\Omega_{\nnu_0}}) \leq J_c(\mathbb
R^2)$ for the cost functional (\ref{cost}).  Otherwise, we divide
$a_{0,0}^\ell$ by $2$ until this requirement is fullfilled.

These values define the prior mean $\nnu_0$ in (\ref{regcost}) and
(\ref{posterior}), which are also the initial guesses shown in
the figures discussed in Sections \ref{sec:blresults} and
\ref{sec:bmcmc}. We select the covariances that render non-positive radii
in (\ref{trigonometric2}) unlikely.
For that purpose, we choose $\Gpr$ as a block diagonal matrix. In
each block, the components decrease with increasing mode number.
Additionally, in our imaging set-up, the uncertainty in the
direction of the incident wave $c_{y,0}^\ell$ is larger than the
uncertainty in $c_{x,0}^\ell$.  
This is due to the increased uncertainty of the topological derivative
in the incident wave, i.e., the $y$-direction.
Specifically, we choose the variances $(\sigma_x^\ell)^2=0.1$ and
$(\sigma_y^\ell)^2=0.2$. 
The variances for the radius mode coefficients $a_{m,0}$ and $b_{m,0}$ are 
inspired by convergence results for deterministic approaches which use 
$H^s$-norms for the radius, that is, weighted $L^2$-norms forcing decay
\cite{hohage2d}.  We set $(\sigma_0^\ell)^2=0.1$ and
$(\sigma_m^\ell)^2 = 0.1/(1+m^2)^s$, $s$ large, $ 1 \leq m \leq M$.
In our tests, we usually set $s=3$. In this way, the prior favors regular
shapes, i.e., shapes with  $r(t)>0$.

Once the prior distribution is defined, we resort to different techniques to 
explore the posterior distribution. Methods for doing this as well as
numerical results are presented in the next two sections. The simplest one
computes a maximum a  posteriori (MAP) estimate and samples the 
linearized posterior distribution to infer properties of the objects which 
generated the data, as illustrated in Fig. \ref{fig2}(c). Green contours
represent probabilities of belonging to a boundary (built from the samples),
whereas the magenta curve represents the MAP point.

\section{Sampling from a Bayesian linearized formulation}
\label{sec:bl}

An approximation of the posterior density (\ref{posterior}), which
builds on tools often available for optimization in deterministic
inverse problems, is the Laplace approximation obtained by
linearization at the maximum a posterior (MAP) point. The approach
first computes the MAP parameter vector $\nnumap$, which
minimizes the negative log likelihood (\ref{regcost}), and then
approximates the posterior distribution by a multivariate Gaussian ${\cal
  N}(\nnumap, \Gpo)$ with posterior convariance matrix $\Gpo =
\mathbf H^{-1}$, where $\mathbf H$ is the Hessian of
(\ref{regcost}) evaluated at $\nnumap$
\cite{Tarantola05,KaipioSomersalo06}. If the parameter-to-observale
map $\mathbf f(\cdot)$ were linear (and with the additive Gaussian
noise and Gaussian prior assumptions made in Section
\ref{sec:bayesian}), this approximation of the posterior distribution
would be exact. In general, the accuracy of this Gaussian posterior
approximation depends on the degree of nonlinearity of $\mathbf f$.

Efficient computation of the MAP point is discussed in
Section~\ref{sec:map}. Once $\nnumap$ is available, the Hessian
can either be computed explicitly or using low-rank approximations for
large or infinite-dimensional parameters
\cite{georglinearized,georgmcmc}.  For a parameter-to-observable map
$\mathbf f(\nnu) = (f(u_{\Omega_{\nnu}}(\mathbf x_j)))_{j=1}^N$ with
measurement operator $f: \mathbb C \rightarrow \mathbb C$, $f(u)=u$
one obtains
\begin{eqnarray} \label{lposterior}
\Gpo = 
\left( {\rm Re}[\mathbf F^{\rm ad} \Gn^{-1} 
\mathbf F] + \Gpr^{-1} \right)^{-1}.
\end{eqnarray}
Here, $\mathbf F$ is the Jacobian matrix of $\mathbf f(\cdot)$ evaluated 
at $\nnumap$ and its adjoint 
$\mathbf F^{\rm ad}= \overline{\mathbf F}^\top$ is the conjugate
transpose of $\mathbf F$. If the measurement operator is 
$f: \mathbb C \rightarrow \mathbb R$, $f(u)=|u|^2$, the first part of
the Hessian, which represents the amount of information learned from the data,
is ${\rm Re}[\mathbf F^{\rm ad} \mathbf M_h \Gn^{-1} \mathbf F]$,
where $\mathbf M_h$ is a real diagonal matrix, see Section
\ref{sec:map}.  Samples from this posterior distribution approximation
${\cal N}(\nnumap, \Gpo)$ can be drawn as
\begin{eqnarray} \label{samplelposterior}
\nnu^{\rm pt} = \nnumap
+ \Gpo^{1/2} \mathbf n,
\end{eqnarray}
where $\mathbf n$ is a vector of independent and
identically distributed (iid) standard normal random values and
$\Gpo^{1/2}$ is a square root of the positive posterior covariance matrix
\cite{georgmarmousi}. 
Analgously, samples from the prior are drawn using 
$\nnu^{\rm pr} = \nnu_0 + 
\Gpr^{1/2} \mathbf n$.

To compute the MAP point, we compared different strategies to
solve the nonlinear least-squares problem (\ref{regcost}): gradient
descent, Newton methods, Gauss-Newton (GN) and Levenberg-Marquardt
(LM) variants.  We detail the procedure finally used next.

\subsection{Computing the MAP point}
\label{sec:map}

Newton methods to minimize a  functional $J(\nnu)$ 
implement the iteration $ \nnu^{k+1} =  \nnu^{k}
- \mathbf H^{-1} \mathbf g$, where $\mathbf H$ is the Hessian and 
$\mathbf g$ is the gradient of $J$. This is equivalent to solving
systems of the form $\mathbf H \boldsymbol \xi^{k+1} = - \mathbf g$,
followed by the update step $\nnu^{k+1}:=\nnu^{k}+\boldsymbol\xi^{k+1}$.
Levenberg-Marquardt approaches  add a `damping' term $ \mu \mathbf I$,
$\mu >0$, to the system matrix. The value of $\mu$ is adjusted at each
iteration. If the objective yields a strong decrease from one
iteration to the next, small values of $\mu$ are used and the method 
resembles a Newton-type scheme. 
When the decrease is slow, large values of $\mu$ are selected and the 
method becomes closer to gradient descent. We use a variant of this 
method proposed by Fletcher \cite{modifiedlm}, which scales this `damping' 
term to allow for larger steps in directions along which the gradient is smaller 
using $ \mu \, {\rm diag}(\mathbf H)$ instead of $\mu\mathbf I$.

For the objective (\ref{regcost}) and the two forms of
the parameter-to-observable map $\mathbf f$, we compute gradients and
Hessians as follows, where for the latter we neglect second-order
derivatives of measurement operators. The resulting Hessian
approximation is sometimes called Gauss-Newton Hessian.  Given a
parameter-to-observable map $\mathbf f(\nnu) =
(f(u_{\Omega_{\nnu}}(\mathbf x_j)))_{j=1}^N$ with measurement operator
$f: \mathbb C \rightarrow \mathbb C$, $f(u)=u,$ we have
\begin{eqnarray*}
\mathbf g(\nnu) := 
{\rm Re} [\mathbf F^{\rm ad}(\nnu)  \Gn^{-1} (\mathbf f(\nnu) - \dn)]
+ \Gpr^{-1} (\nnu
-\nnu_0), \\
\mathbf H^{\rm GN}(\nnu)  :=  {\rm Re}[\mathbf F^{\rm ad}(\nnu) \Gn^{-1} \mathbf F(\nnu)] + \Gpr^{-1},
\end{eqnarray*}
where here and in the following, $\mathbf F$ denotes the Fr\'echet derivative  
of the parameter-to-observable map
$\nnu \rightarrow (u_{\Omega_{\nnu}}(\mathbf x_j))_{ j=1}^N,$
$u_{\Omega_{\nnu}}$ being the solution of (\ref{forward}).
The characterization of this operator in terms of the solutions
of auxiliary boundary value problems, as well as the way to calculate
$\mathbf F,$ are discussed in \ref{sec:chfrechet}. 
For a measurement operator $f: \mathbb C \rightarrow \mathbb R$, 
$f(u)=|u|^2,$ and diagonal constant $\Gn$, gradient and Hessian are
\begin{eqnarray*}
\mathbf g(\nnu) :=  {1 \over \sigma_{\rm noise}^2}
{\rm Re} [\mathbf F^{\rm ad}(\nnu)  \mathbf M_g (\mathbf f(\nnu) - \dn)]
+ \Gpr^{-1} (\nnu
-\nnu_0), \\
\mathbf H^{\rm GN}(\nnu)  :=   {1\over \sigma_{\rm noise}^2}
{\rm Re}[\mathbf F^{\rm ad}(\nnu)  \mathbf M_h
\mathbf F(\nnu) ] +  \Gpr^{-1},
\end{eqnarray*}
where $\mathbf M_g$ and $\mathbf M_h$ are diagonal matrices defined as follows: 
\begin{equation*}
\mathbf M_g= {\rm diag}\left[(2 u_{\Omega_{\nnu}}(\mathbf x_j))_{j=1}^N \right], \qquad
\mathbf M_h = {\rm diag}\left[(6|u_{\Omega_{\nnu}}(\mathbf x_j)|^2
- 2 \dnj)_{j=1}^N \right].
\end{equation*}

To summarize, we use the Levenberg-Marquardt method with scaled
diagonal matrix and the Gauss-Newton Hessian to compute the MAP point.
That is, starting from $\nnu^0 = \nnu_0$, we use the iteration $
\nnu^{k+1} = \nnu^{k} + \boldsymbol \xi^{k+1}$, where 
$\boldsymbol \xi^{k+1}$ is the solution of
\begin{eqnarray*}
\left(\mathbf H^{\rm GN}_{\lambda_k}(\nnu^k) + \mu_k 
{\rm diag}(\mathbf H^{\rm GN}_{\lambda_k}(\nnu^k) ) \right)
\boldsymbol \xi^{k+1} = - \mathbf g_{\lambda_k}(\nnu^k).
\end{eqnarray*}
Here,  the subscript $\lambda_k$  indicates that we multiply $\Gpr^{-1}$
by a factor $\lambda_k$ in the initial steps to balance the different orders of 
magnitude of the two terms defining the cost $J$ in
(\ref{regcost}). The initial 
value $\lambda_0 =0.1 \sigma_{\rm noise}^{-2}$  decreases each iteration
by a factor $2/3$ until it reaches the value $\lambda=1$ corresponding to
(\ref{regcost}). We set $\mu_k = 10^{-3}$ for each iteration, and check
if the  functional decreases sufficiently. If it does not, this
value is increased by a fixed factor.
This usually only happens for a few steps, and $\mu$ is 
bounded from above by $1$. The iteration stops when the
difference between the new value of the cost and the
previous one is smaller than $\tau  \sigma_{\rm noise}^{-2}$, where
$\tau$ is a tolerance (usually $\tau = 10^{-5}$).
The final parametrization $\nnumap$ is considered the MAP
point. Computing this MAP point typically requires 15-25 iterations
for problems with single objects.

Notice that a star-shaped parameterization of an object is uniquely
defined when its center is fixed.  However, star-shaped objects
(circles or ellipses, for instance) may be parameterized with
different centers, at least when infinitely many coefficients are used
in the expansion.
If we consider only the cost
functional (\ref{cost}), we may encounter different minima defining
the same or very similar objects.  Some optimization strategies
\cite{jcp19} overcome this issue by enforcing that the center of the
object must be the center of mass. In our Bayesian approach, the regularization
(which results from our choice of prior) in (\ref{regcost}) results in a
preference for one parametrizations for the same object.

\subsection{Object detection with quantified uncertainty}
\label{sec:blresults}

Once we have computed the MAP point, we can apply the
Bayesian linearized framework to object detection. We have
performed our tests in the set-up shown in Fig.\ \ref{fig1}(a), adapted 
to either acoustic or light holography; these settings are described in
detail in \ref{sec:physical}. 
We place detectors at $\mathbf x_j=-5  +0.05  j$, $j=0,\ldots,200$,
on a uniform grid of step size $0.05$.
The subgrid with step size $0.1$ is used to compute the MAP point, whereas 
the data at intermediate detectors is used to generate the topological priors.
The same splitting is used for our MCMC sampling experiments in the
next section.

\begin{figure}[h!]
\centering
\includegraphics[width=11cm]{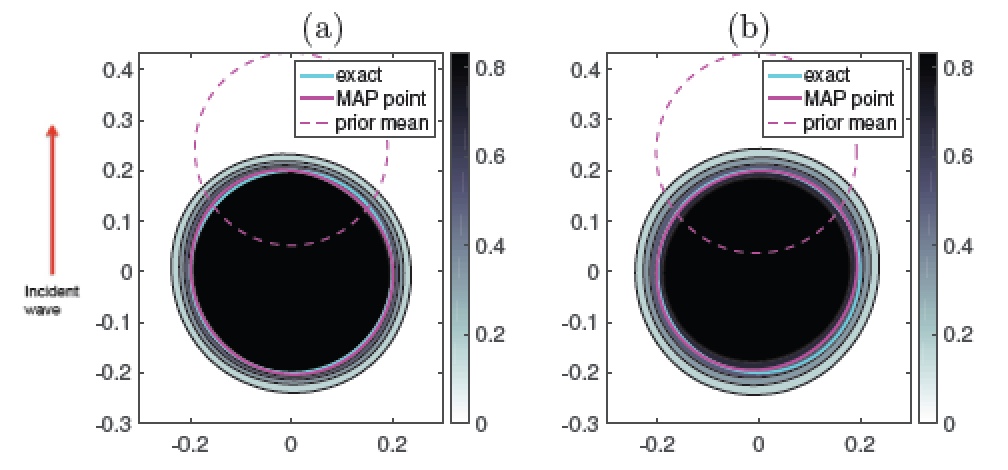} 
\caption{Linearized Bayesian solution for synthetic measurements of
  complex amplitudes using a sphere centered at $(0,0)$ with radius
  $0.2$ with noise (a) $1 \%$ and (b) $5 \%$.  Cyan curves show the
  true object used to generate the synthetic data, whereas the
  obtained MAP points are shown in magenta. Dashed curves are the
  prior means constructed by topological methods, which are also used
  as initial guesses in the optimization.  Green contours show the
  probability of points to belong to the object.  Parameter values are
  $\kappa_{\rm i}=15.12$, $\kappa_{\rm e}=12.56$, $\beta=1$, and
  $M=5$.  The arrow indicates the incidence direction. Detectors are
  placed at a distance $5$ in that direction, as in
  Fig.\ \ref{fig1}(a).  }
\label{fig3}
\end{figure}

\begin{figure}[h!]
\centering
\includegraphics[width=14cm]{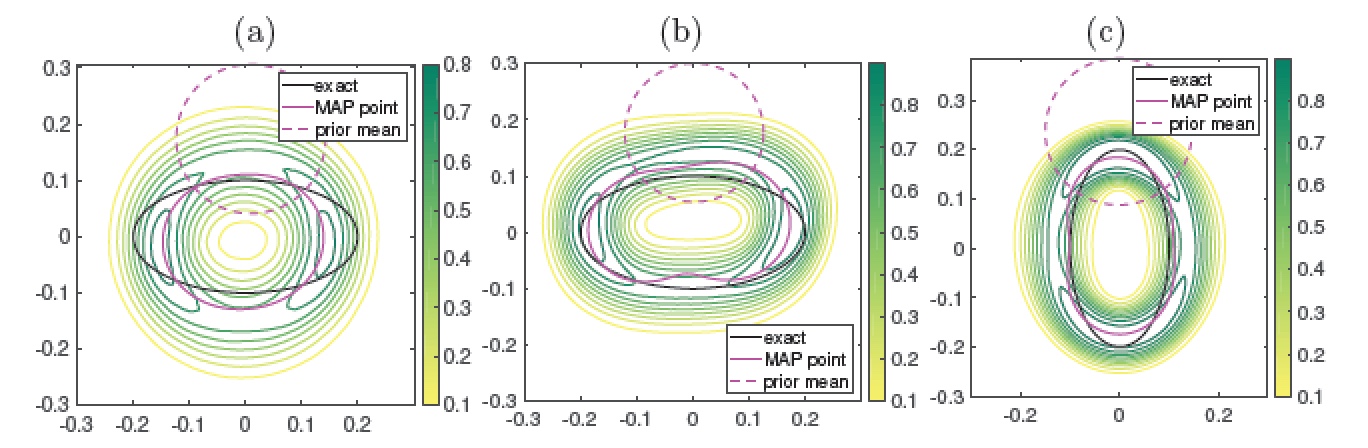}
\caption{(a)-(b) Inferred shapes for an ellipse centered at $(0,0)$
  with semi-axes of lengths $0.2$ and $0.1$ in the $x$ and $y$
  directions, respectively, when (a) $\kappa_{\rm i}=15.12$,
  $\kappa_{\rm e}=12.56$, and (b) $\kappa_{\rm i}=25.12$, $\kappa_{\rm
    e}=20.56$, with $\beta=1$, $M=5$.  (c) displays results for the
  same parameters as in (b) orienting the ellipse in the $y$
  direction.  Green countours show the pointwise marginal for each
  point on the curve.  Data: Complex amplitudes at detectors placed at
  a distance $5$ in the incidence direction, with $5 \%$ noise.  }
\label{fig4}
\end{figure}

\begin{figure}[h!]
\centering
\includegraphics[width=11cm]{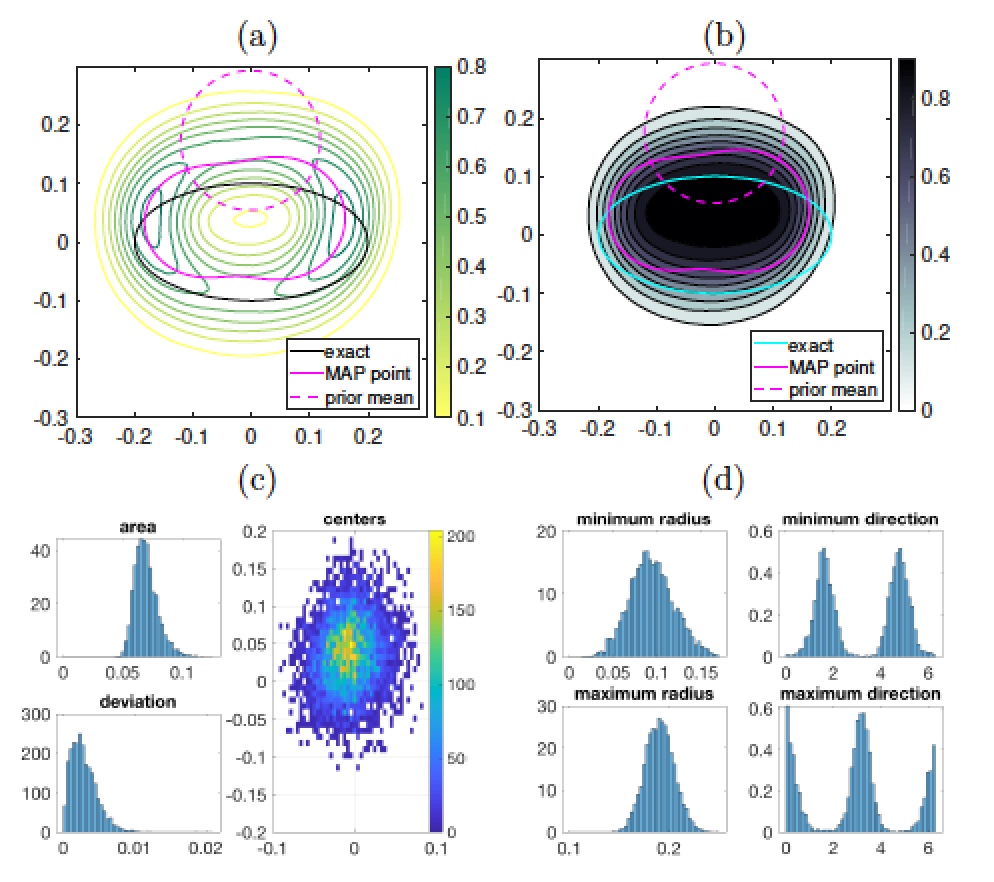} 
\caption{Inferred shapes from complex data with $10 \%$ noise for the
  physical parameters $\kappa_{\rm i}=25.12$, $\kappa_{\rm e}=20.56$,
  $\beta=1$, and $M=5$.  The true shape is an ellipse centered at
  $(0,0)$ with semi-axes of length $0.2$ and $0.1$ in the $x$ and $y$
  directions, respectively.  Contours represent the contours for the
  pointwise marginal for points on the curve in (a) or the probability
  of being inside the object in (b). The histograms in (c) and (d) are
  discrete approximations of the densities for the distribution of the
  area, the deviation, the center of mass, the minimal and maximal
  radius and the angle of the direction of the minimal and maximal
  extension of the inferred objects.}
\label{fig5}
\end{figure}

Figures \ref{fig2}-\ref{fig5} show the results for two sets of wavenumbers, 
namely $\kappa_{\rm i}=15.12$, $\kappa_{\rm e}=12.56$, and 
$\kappa_{\rm i}=25.12$, $\kappa_{\rm e}=20.56$.
The incident wave is $u_{\rm inc}= e^{\imath \kappa_{\rm e} y}$, where the
$y$-coordinate is oriented as in Fig.\  \ref{fig1}(a).

The topological approach described in Section \ref{sec:priors} is used
to define the priors, and the prior mean is used as initialization to
compute the MAP point. The functional (\ref{regcost}) is minimized as
explained in Section \ref{sec:map} to obtain the MAP
point. Linearizing about it and approximating the posterior
probability by a Gaussian, we generate samples by means of
(\ref{samplelposterior}). Based on them, we compute probabilities for
points to lie inside the object (e.g., Fig.\ \ref{fig4}), the
pointwise marginals of the countour (e.g., Fig.\ \ref{fig5}), as well
as the marginal distributions of the centers of mass, the maximum and
minimum distances from the center to the object border and their
orientation. A small number of samples have negative radius as defined
by (\ref{trigonometric1})-(\ref{trigonometric2}) resulting in shapes
with loops, which we discard for these calculations.

Figure \ref{fig3} compares the results for a circle-shaped object as the 
noise magnitude increases from $1\%$ to $5\%$, for the lower wavenumbers.
The position, size and shape are inferred with small uncertainty.
Switching to ellipsoidal shapes, the values
of $\kappa_{\rm i}$, $\kappa_{\rm e}$ need to be increased for a more precise description 
of the dimensions and the orientation of the object, see Figure 
\ref{fig4}. However, smaller values locate the true center more accurately.
A similar phenomenon is observed in 3D \cite{jcp19}.
Improved resolution has also been achieved increasing 
the wavenumber in acoustic 2D settings \cite{greengard} with full 
aperture far field measurements, however,  the number
of incident directions and sampling points is increased too. Here,
we keep the same single incident direction and limited aperture 
sampling points.
The inferred shapes are still reasonable for higher noise
magnitudes. Figure \ref{fig5} illustrates the results for an
ellipse, including statistics for its center of mass, size and
orientation, with $10\%$ noise.  
Similar results are obtained when halfing the number
of detectors (these results are not shown here). If we reduce the number of detectors
further, the distance between them becomes larger than
the object size. As a consequence, the MAP point looses the elliptical shape, but
its location and size are still adequate.
In Figure \ref{fig5}, we additionally
study the uncertainty in the inferred objects
in terms of quantities that do not depend on the parametrization, namely the
area, the center of mass and the minimal and maximal radii and their
directions. The area for a closed curve ${\cal C}$ is computed as
$\int_{\cal C} r(s) ds$. The center of mass of a curve with
the parameterization $\mathbf q(t) = (x(t),y(t))$,
$t \in I$, is  ${\int_I \mathbf q(t) |\mathbf q'(t)| dt/ \int_I 
|\mathbf q'(t)| dt }$, where $|\mathbf q'(t)| dt = \sqrt{x'(t)^2 + y'(t)^2}  
dt = ds$ is the differential of arch length.
The deviation from a circular shape is $\int_{\cal C} |r(s) - r_{av}| ds$, where
$r(s)$ is the distance to the center of mass and $r_{av}$
the average radius. The maximum and minimum radii are the longest and 
shortest distances of the curve to the center of mass, respectively. The angle 
between the $x$-axis and the direction of longest and shortest radius is denoted
as maximum and minimum direction in Figure \ref{fig5}.
Figures \ref{fig3}-\ref{fig5} use as data complex amplitudes measured 
at $101$ detectors placed at a distance $5$ of the object in the 
incidence direction. As commented earlier, Figure \ref{fig2}(c) 
uses real intensities  at the detectors as 
measured data. Unlike ellipsoidal shapes, these three objects 
can be exactly parametrized for a small number
$M=5$ of modes to represent the object geometry. The reconstruction of 
the configuration is  reasonable. 

The computational cost of this procedure is moderate. The MAP points
are usually obtained in about $15$-$25$ iterations, and sampling with
the expression (\ref{samplelposterior}) requires no additional forward
solves besides computing the Hessian at the MAP point.  The figures
discussed here use $10,000$ samples.  However, we have introduced some
approximations in the Bayesian inference process: first, when
approximating the posterior by a Gaussian, second, when approximating
the full Hessian with the Gauss-Newton Hessian at the MAP point.  To
assess the validity of the procedure, we  
will compare to results obtained by MCMC sampling of the whole posterior 
distribution.  

\section{Sampling with Markov Chain Monte Carlo}
\label{sec:bmcmc}

Markov Chain Monte Carlo (MCMC) techniques have the potential to fully
explore and statistically characterize the posterior distributions without linearizing the
parameter-to-observable map.
A Markov chain is a sequential stochastic process, which moves from one state to another 
within an allowed set of states: $X^{0} \longrightarrow  X^{1} \ldots \longrightarrow  
X^k \ldots$. To define a Markov chain we need three elements:
1) the state space, that is, the set of states $X$ the chain is allowed to reach,
2) the transition operator $p(X^{k+1} | X^k)$ which establishes the probability 
of transitioning from state $X^k$ to $X^{k+1}$, and
3) the initial distribution $\pi_0$ which defines the initial probability of being 
in any of the possible states.
To generate a Markov chain, one
moves from one state to another guided by the transition operator 
$p(X^{k+1} | X^{k})$.

In our context, we wish to sample posterior distributions $\pi$ by 
means of Markov chains, and a natural choice for $\pi_0$ is the prior
distribution.
There are many MCMC variants for sampling posterior distributions
\cite{optimizationbayes}, which are based on different transition
operators.
MCMC algorithms often suffer from ``slow mixing'', i.e., the underlying
Markov Chain takes (too) many samples to explore the parameter space
and thus provide a good characterization of the distribution we aim at
sampling. For large scale problems near a continuous limit, preconditioned 
Crank-Nicholson and stochastic Newton variants have been shown successful
\cite{dunlopthesis,georgmcmc, CotterRobertsStuartEtAl13}. 
However, when the target distribution is multimodal, these samplers may fail to
jump from one mode to another.
Thus, we resort to affine invariant MCMC 
samplers working with several chains \cite{goodmanweare}
because their strategy to generate new proposals reduces the occurrence 
of samples with negative radii during the sampling process and
the use of many chains allows us to handle multimodal distributions.

In our imaging set-up, the initial distribution $\pi_0$ is the
truncated prior distribution, where we neglect normalization factors
\begin{eqnarray} \label{mcmcprior}
\pi_{0}(\nnu) = \left \{
\begin{array}{ll}
0, & \mbox{if  intersections}, \\
\exp(- {1 \over 2} \left(\nnu - \nnu_0\right)^t 
\Gpr^{-1} (\nnu - \nnu_0) ), & \mbox{otherwise}.
\end{array} \right.
\end{eqnarray}
By intersection we mean that the radius of the star-shaped curves
with coefficients $\nnu$ vanishes at at least one point and thus the
curve is degenerate, and, for multiple objects, that we have 
intersecting or nested boundaries.
Disregarding  the normalization constants, which do not play a role in
MCMC methods, the posterior to be sampled is
\begin{eqnarray*}  
\pi(\nnu) = \pi_{0}(\nnu) \exp\big( -{1\over 2} 
(\overline{\dn - \mathbf f(\nnu)})^t 
\Gn^{-1} 
(\dn - \mathbf f(\nnu))
\big).
\end{eqnarray*}
For the algorithm to sample the posterior distribution we follow
\cite{goodmanweare}, which we summarize for completeness:
\begin{itemize}
\item Initialization: 
Generate the initial positions of the walkers $X_0^{w} \in \mathbb R^d$, 
$w=1,\ldots,W$ by sampling from the  prior distribution $\pi_0$.
Choose the acceptance parameter $a$ (typically
$2$) and the number of samples $K$.
\item for each step $k=0,\ldots,K-1$, evolve the walkers
 $w=1,\ldots,W$ as follows:
  \begin{itemize}
   \item  Draw a walker $X^q_k$ at random from the set of walkers
   $\{ X^j_k \}_{j \neq w}$.
   \item Choose a random value $z_w$ from the distribution $g(z)= {
   1\over \sqrt{z}}$  when $z \in [1/a,a]$, zero otherwise.
   \item Calculate proposition $X_{\rm prop}^w = X_k^q  + z_w
   (X_ k^{w} - X_k^{q})$.
   \item Calculate $s = z_w^{d-1} {\pi(X_{\rm prop}^w) \over \pi(X_ k^{w})}.$
 Calculate $s = {\rm Min} \left( 1,s \right).$
   \item Draw $r$ with probability ${\cal U}(0,1)$. If $r \leq s$ set $X_{k+1}^w = X_{\rm prop}^w,$  
   otherwise set $X_{k+1}^w = X_k^w.$
 \end{itemize}
\item  Final result: The Markov chains $\{ X_0^w,\ldots, X_K^w \}$,
for all the walkers $w=1,\ldots,W$.
\end{itemize}
The number of walkers is chosen as at least $W>2d$ to increase mixing 
and to take advantage of parallelization \cite{hammer}.

\begin{figure}[h!]
\centering
\includegraphics[width=14cm]{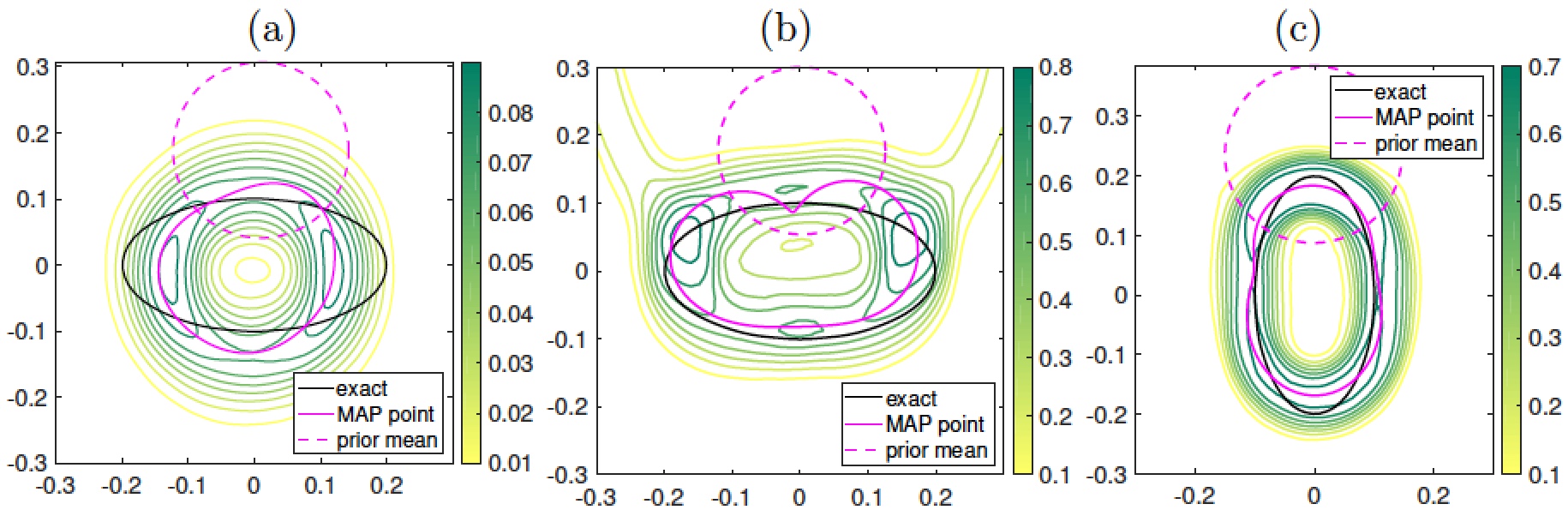} 
\caption{Same as Figure \ref{fig4} but with MCMC sampling. Sampling
  parameters are $W=200$, $K=200$, $B=35000$ and $a=2$ in (a), and
  increased to $K=500$ and $B=55000$ in (b)-(c). }
\label{fig6}
\end{figure}

\begin{figure}[h!]
\centering
\includegraphics[width=14cm]{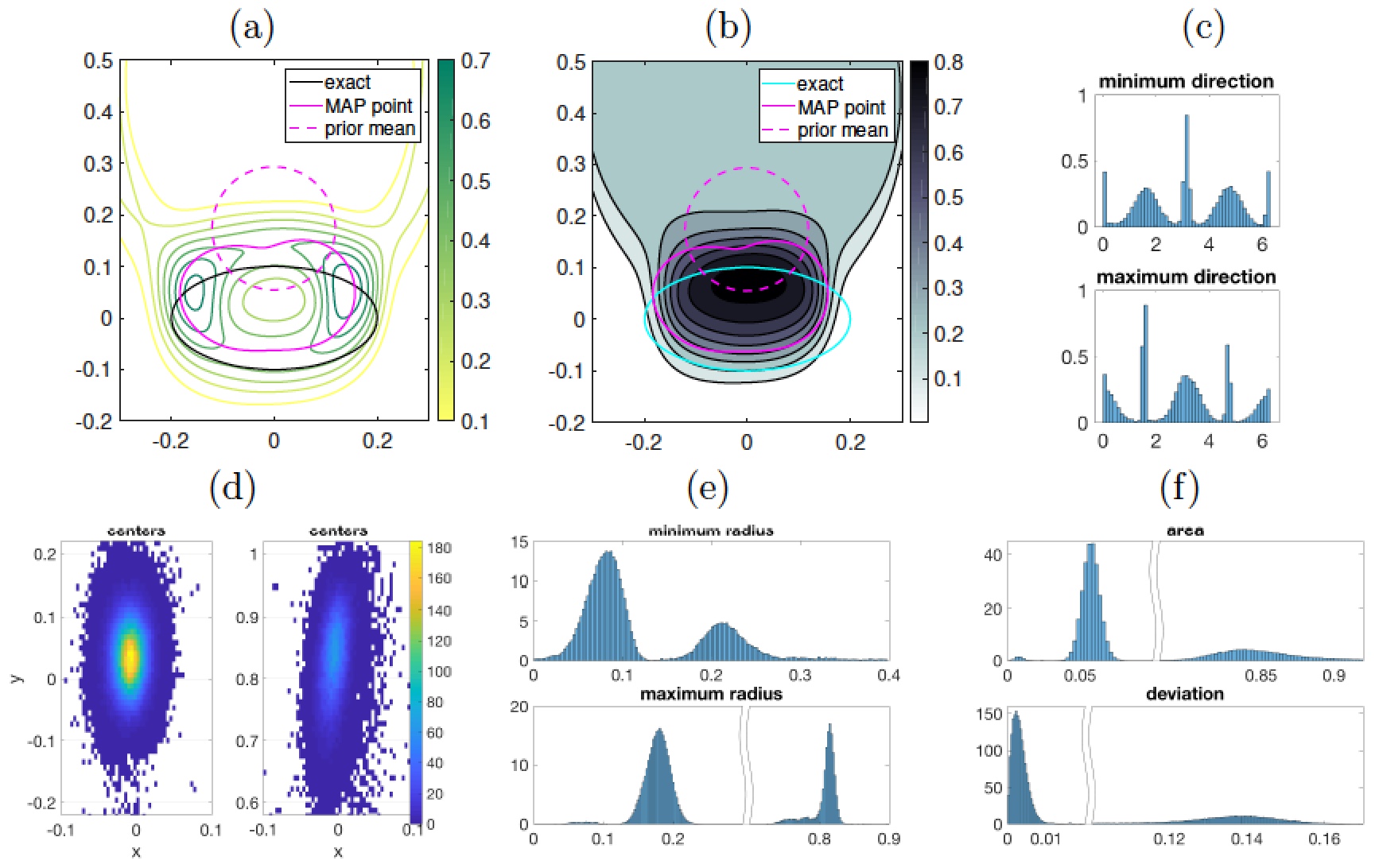}

\caption{Same as Figure \ref{fig5} but with MCMC sampling. Sampling
  parameters are $a=2$, $W=800$, $K=1500$, $B=80000$. Due to the
  bi-model nature of the posterior distribution, only parts of the
  figures are shown. See Figure \ref{fig9} for a study explaining of
  the bimodality of the posterior.}
\label{fig7}
\end{figure}

With this  algorithm, we revisit the configurations studied by 
Bayesian linearized techniques. Figures \ref{fig6}-\ref{fig7} are the 
counterparts  of figures \ref{fig4}-\ref{fig5}, generated by MCMC sampling. 
We have used Gelman-Rubin tests \cite{gelman} to check that the
sample distributions  under study are properly converged.
Notice that  in both cases we have rejected samples with negative radius. 
The results are similar but the computational cost is
higher, several tens of thousand samples have been 
used, depending on the shapes and wavenumbers. Each sample
requires solving a forward problem to evaluate the likelihood.
Magenta curves in this case represent the  
sample with maximum probability, which is an 
approximation to the MAP point. The approximation improves
increasing the number of samples.

\begin{figure}[h!]
\centering
\includegraphics[width=14cm]{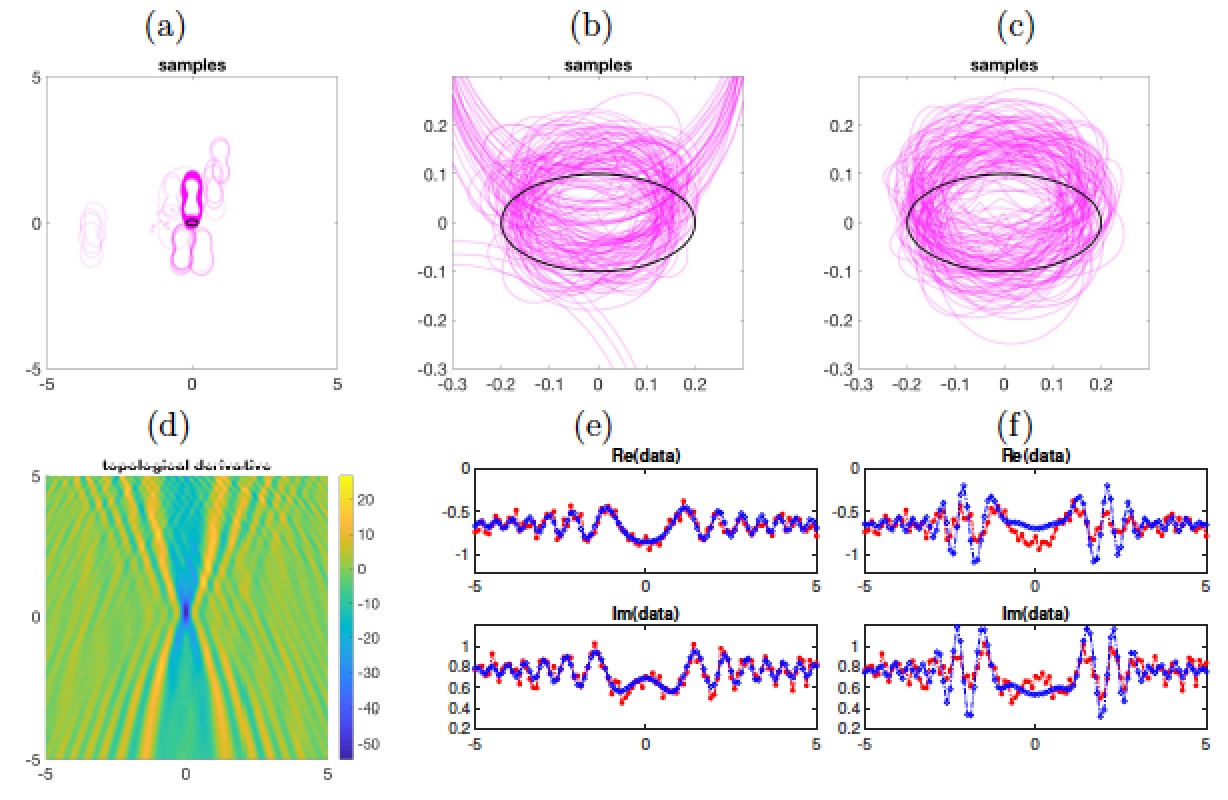} 
\caption{Samples from the problems discussed in Figures \ref{fig5} and
  \ref{fig7}. Shown in (a) are the last 800 walker states generated by
  the MCMC algorithm in the imaging region $[-5,5]\times [-5,5]$, in
  (b) zoom-ins showing $100$ of such samples around the true object,
  and in (c) $100$ samples generated by Bayesian linearized
  sampling. Using the linearized approach, all samples wrap around the
  true object (black curve), while samples from the posterior obtained
  using MCMC show that other shapes are consistent with the data. To
  study this further, in (d) the topological derivative of the cost
  functional (\ref{cost}) is shown, and in (e)-(f) a comparison of the
  data (asterisks) with observations generated from a sample wrapped
  around the object (circles in (e)) and by one of the large,
  elongated samples (circles in (f)).}
\label{fig9}
\end{figure}

Comparing Figure \ref{fig6} and Figure \ref{fig4}, we see that 
the approximations to the MAP point are similar. However, there
are differences in the lower probability contours around the MAP point
for the ellipse oriented orthogonally to the incidence direction and
the larger values of $\kappa$, see panel (b). These low probability
features are converged as when increasing the number of Monte Carlo
samples, they remain mostly unchanged.
We analyze this phenomenon in more detail in Figures
\ref{fig7}-\ref{fig9} for larger noise. Comparing Figure \ref{fig7} to
Figure \ref{fig5}, the MAP and expected curves show again reasonable
agreement. However, the Bayesian posterior becomes bimodal, as seen in
the histograms for centers, radii, orientation and area.  These
observations persist varying the number of walkers $W=200, 400, 800$,
steps $S=500, 1000, 3000, 7000, 14000$, and acceptance rates $a=2,
1.1, 2.5$. 

To further study this, Figure \ref{fig9}(a) displays a collection of
samples. Most of them wrap around the true object, as shown in panel
(b). However, a significant number of samples is oriented in the
incidence direction of the waves, along the $y$ axis, orthogonal to
the true orientation of the object.  Tracking the evolution of the
initial walkers, we observe that when they are above a certain size,
the proposed curves evolve toward this second family, whereas smaller
shapes get closer to the true object.  Comparing Figure \ref{fig9}(d)
with Figure \ref{fig9}(a), we notice that samples concentrate in
regions of large negative values for the topological derivative. Some
concentrate around the spot where largest negative values are
attained, marking the true object location, while the other elongate
along stripes of less negative values. As commented in Section
\ref{sec:topologicalfields}, we expect the physical cost functional
(\ref{cost}) to decrease by removing regions of negative values of the
topological derivative, that is, placing objects in them.  Panels
(d)-(e) compare the synthetic data used in the simulations to the
measurements corresponding to samples wrapped around the true object
and to elongated samples.

Elongation and loss of axial resolution are aberrations present in
traditional holographic reconstruction techniques based on numerical
backpropagation \cite{backpropagation} due to the use of only one
incident direction.  They indicate the potential presence of
additional local minima in the original cost functional and highlight
ambiguity due to ill-posedness, in particular if the data contains
larger noise levels.

\begin{figure} 
\centering
\includegraphics[width=11cm]{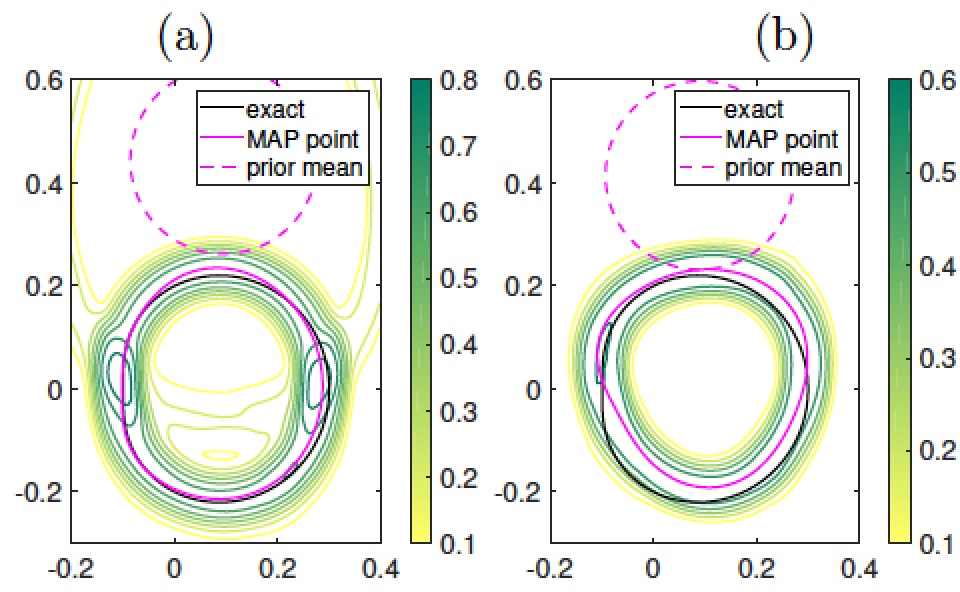}
\caption{Results for an egg-like object inferred from
complex data $u(\mathbf x_j)$ (a), and from the modulus 
 data  $|u(\mathbf x_j)|^2$ (b), $j=1,\ldots,101$, with $5\%$ noise.
Parameters are $\kappa_{\rm i}=24.79$, $\kappa_{\rm e}=20.6$,
$\beta=1$ and $M=5$. Sampling parameters 
are $a=2$, $W=200$ and $K=3000$.
}
\label{fig10}
\end{figure}

Figure \ref{fig10} compares the results for a asymmetric
egg-like shapes when we replace complex valued data by just
intensities.  As can be seen, when only using intensities, the
uncertainty in the shape increases.

\section{Varying the material properties}
\label{sec:material}

\begin{figure}[h!]
\centering
\includegraphics[width=11cm]{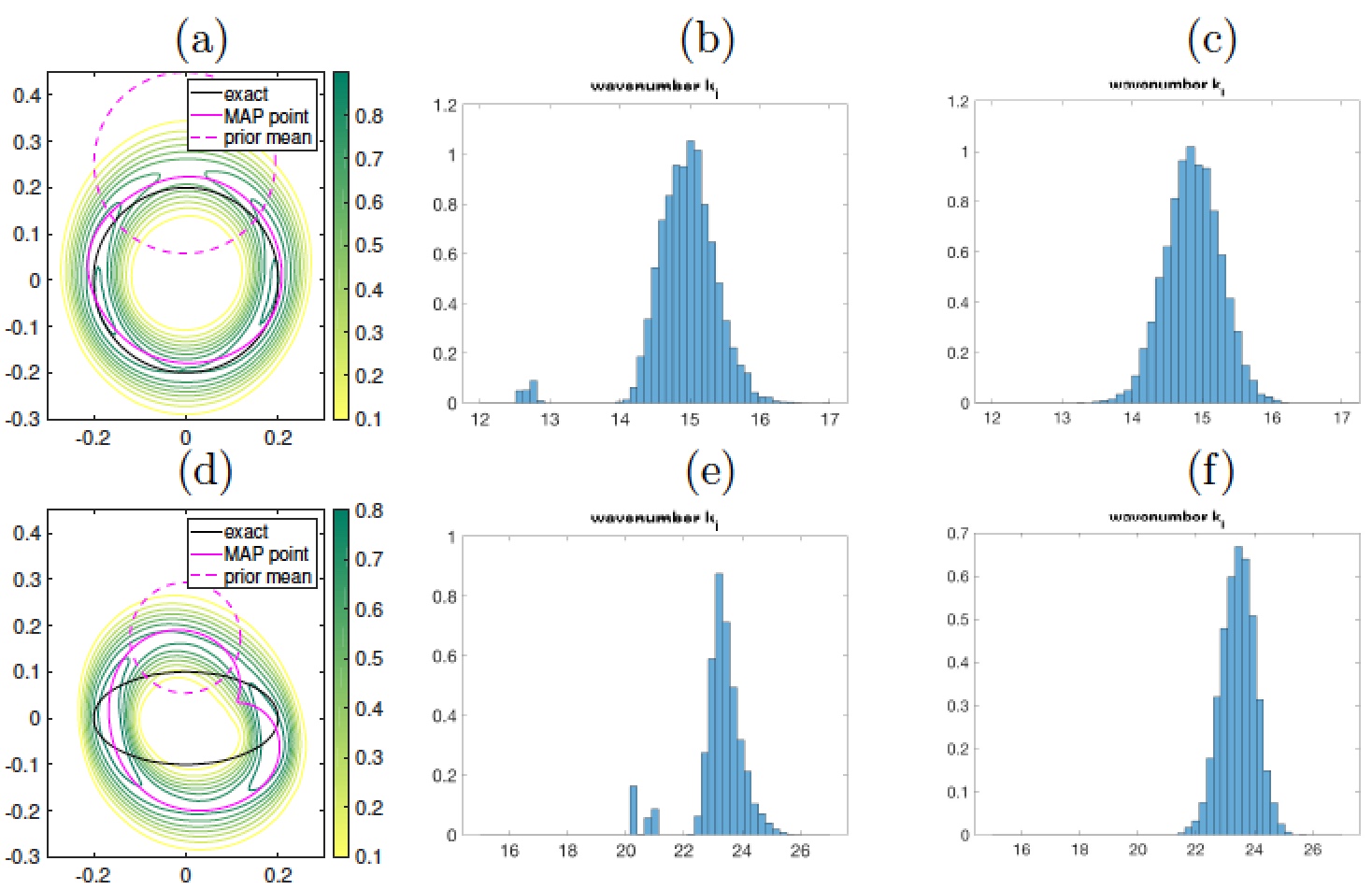}
\caption{(a) Same as Figure \ref{fig3}(b) but with MCMC sampling for
  unknown $\kappa_{\rm i}$, shown in (b). Shown in (d) is the
  counterpart of Figure \ref{fig4}(b) using intensity data with MCMC
  sampling for unknown $\kappa_{\rm i}$, shown in (e).  (c) and (f)
  are the equivalent of (b) and (e) obtained by linearizing about the
  MAP point.  We use the same MCMC sampling parameters as for the
  previous figures. The  true values for $\kappa_{\rm i}$ are
  $\kappa_{\rm i}=15.12$ in (a)-(c) and $\kappa_{\rm i}=25.12$
  in (d)-(f). Expected values are $\kappa_{\rm i}=15.04$ and 
  $\kappa_{\rm i}=23.18$, respectively.
  }
\label{fig12}
\end{figure}

In the previous sections, we fix the object properties and aim at inferring
their geometry. In practice, one may have to infer also material
constants entering the equations such as $\kappa_{\rm i}$. In the MCMC
framework, this can be handled by adding an additional element to the
prior. The forward problem has one parameter in addition to those
defining the shapes, which enters the equation.  Hence, we sample
with respect to one more parameter. The results are shown in
Fig.\ \ref{fig12}.

In the Bayesian linearized framework, we also include a new parameter
in the prior. When optimizing, we simply compute derivatives with
respect to one parameter more. The Fr\'echet derivative with respect
to $\kappa_{\rm i}$ is then the solution of a boundary value problem
obtaining differentiating the Helmholtz equations with respect to
$\kappa_{\rm i}$, see \ref{sec:coefficients}. Once this is done, a
similar iteration to that proposed in Section \ref{sec:map} provides a
MAP estimate. Sampling by means of (\ref{samplelposterior}) we obtain
Fig.\ \ref{fig12}(c) and (f).

\section{Model selection when the number of objects is unknown}
\label{sec:model}

In the previous sections, we inferred object shapes from observational
data assuming that the number of objects is known. Here, we are
considering the number of objects $L$ as an additional variable to be
inferred. The number of objects is substantially different from other
variables that characterize the object geometry because it is discrete
and changing it amounts to the addition or removal of blocks of
parameters defining object components.  Thus, the number of objects
defines different models with different numbers of parameters
\cite{modelselection}. The probability of the data given the
model, that is, the number of objects $m$, is computed by
integrating over the model parameters:
\begin{eqnarray}
p(\dn | m) = \int_{\nnu} p(\dn | \nnu) p(\nnu | m) d \nnu,
\label{evidence}
\end{eqnarray}
where $p(\dn | \nnu)$ is the likelihood (\ref{likelihood}) and $p(\nnu
| m)$ the prior (\ref{mcmcprior}) given the model.  It can be
evaluated by sampling one of them, $p(\nnu | m)$ for instance, and
computing a Monte Carlo estimate based on these samples $\nnu_i$,
$i=1,\ldots,S$, amounting to ${1\over S}\sum_{i=1}^S p(\dn | \nnu_i)$,
\cite{optimizationbayes}. The resulting quantity is called the
evidence for model $m$ \cite{penny06}.
We have implemented this procedure for the synthetic observations
coming from the configuration with $3$
objects shown in Figure \ref{fig2}. Computing the evidence for
models with $m=1,2,3,4$, we find a clear
maximum for $m=3$,
i.e., the true number of objects. Once we have selected
a number of objects, we sample the posterior distribution (\ref{posterior})
using the MCMC techniques described in Section \ref{sec:bmcmc} to
infer the expected objects, which resemble the MAP points in Figure \ref{fig2}(c).

\section{Summary and conclusions}
\label{sec:conclusions}

We have developed a Bayesian framework for object detection
which uses topological methods to generate priors. 
In this approach, objects are represented by star-shaped
parameterizations.
Assuming the number of objects and their material properties are known,
we compute the `maximum a posteriori' (MAP) estimate
for the parameters defining centers and radii through minimization of
the proper cost functional. Linearizing the parameter-to-observable
map about the MAP point, one can generate
samples of the Laplace approximation to the posterior distribution
to quantify the uncertainty in the object location
and its shape at a low computational cost.  We test the scheme in
2D holography imaging set-ups, for both acoustic and light waves.
The former uses complex fields measured at detectors as
data, whereas in the latter only real intensities are available.
In these set-ups, wave fields are governed by transmission
boundary problems for Helmholtz equations and the incident waves 
reduce to a single beam. For small noise magnitudes,
many shapes can be inferred with moderate uncertainty for a wide range of
wavenumbers. As the magnitude of the noise increases,
larger wavenumbers (i.e., smaller frequencies)
provide better shape descriptions whereas smaller wavenumbers
approximate their location more accurately.
Even for noise magnitudes of $10\%$ and higher,
topological derivatives are useful to generate priors, and the
Bayesian approach with complex data results in
parameter inference with moderate uncertainties.
This is particularly relevant for 
acoustic imaging set-ups, in which the magnitude of the noise is
larger.

To assess the performance of Bayesian linearized methods, i.e.,
Laplace posterior approximations, we
compared the results with those obtained by direct MCMC sampling 
of the full posterior distribution, finding reasonable agreement.
However, MCMC
sampling provides more complete insight into the structure of the posterior
distribution, which become  multimodal for larger noise levels.
Comparisons are made under the assumption of piecewise constant
material properties, which allows us to use fast boundary element 
schemes to solve forward and Fr\'echet problems. Also, working with 
2D  star-shaped objects represented in a trigonometric basis, we only
need to infer  $14-50$ parameters for sets of $1$-$3$ objects. 
When the number of objects is also unknown, we used
a Bayesian model selection strategy to obtain information on the 
number of objects, as well as on their location and shapes. 

Our methods extend to 3D set-ups, involving a larger number of
parameters.  Whereas Bayesian linearized techniques may yet be
efficient once the derivatives needed to implement optimization
techniques are characterized as solutions of specific boundary
problems, the computational cost of direct MCMC sampling
increases. Here, computing a MAP point requires 15-25 Newton-type
descent steps.  Each iteration implies solving a small linear system,
as well as boundary value problems for the derivatives. Instead, MCMC
requires solving a large amount of boundary value problems.  The
problem becomes even more difficult if we allow for spatial variations in the material
parameters.  Tempering approaches \cite{tombayesian} may help
  to reduce that cost.

Topological sensitivity is a flexible tool to obtain
prior information. It uses an explicit formula which
accommodates limited aperture data,  only one incident wave, large noise, 
and measurements of the full complex field, or functions of it, such 
as intensities. Moreover, it has an interpretation
as a measure of the decrease of the physical cost functional which 
provides insight into  
the  multimodality of the posterior. 
However, the proposed Bayesian framework is independent of the
choice of the prior.
Other initialization procedures, such as linear or direct sampling
\cite{cakonilsm,dsm}, may be used to construct priors when the 
combination of receivers and incident directions provides 
complex amplitude data widely distributed around the object.

Here, we mostly consider tests in which the material properties 
are considered known. Assuming material properties to be 
characterized by constant parameters, the proposed framework
can be adopted by including a few additional parameters. 
To consider more general spatially variable material properties,
one could combine these methods with those developed in
 \cite{georglinearized,georgmarmousi} and implement coupled
BEM-FEM or spectral-FEM solvers as in \cite{sims18,jcp19}. 
The methods would extend to wave fields governed by systems
different from Helmholtz equations, provided characterizations for the
derivatives and adequate solvers for the boundary value problems
involved are available.  
Finally, note that the objects in our tests are assumed to be
stationary. Time-dependent Bayesian methods to track moving 
contours are discussed in \cite{tannenbaum}.

\ack A.\ Carpio and S.\ Iakunin acknowledge partial support from the
FEDER/MICINN - AEI grant   MTM2017-84446-C2-1-R. 
G.\ Stadler acknowledges partial support from
the US National Science Foundation grant \#1723211 and from KAUST
under Award \#OSR-2018-CARF-3666.
A.\ Carpio thanks R.E.\ Caflisch for hospitality during a
sabbatical stay at the Courant Institute, NYU, and D.G.\ Grier for an
introduction to acoustic holography.

\appendix
\section{Characterization of Fr\'echet derivatives}
\label{sec:chfrechet}

In this Appendix, we study Fr\'echet derivatives for transmission problems
of the form (\ref{forward}), where the total wavefield $u=u_{\rm inc} + u_{\rm sc}$ in 
the exterior region  $\Omega_{\rm e}=\mathbb{R}^{2}\setminus \overline{\Omega}_{i}$ 
and  the transmitted wave field  $u=u_{\rm tr}$ in a bounded smooth inclusion
$\Omega_{\rm i}$.
We require Fr\'echet derivatives of the solutions with respect to the
parameters defining $\Omega_{\rm i}$ and with respect to $\kappa_{\rm
  i}$. Fr\'echet derivatives with respect to the domain are calculated
using integral equations in \cite{hohage2d,hettlichfrechet}. We recall
these expressions here and give a different proof by variational
methods. Similar, but simpler, variational arguments yield expressions
for the Fr\'echet derivatives with respect to problem coefficients,
for similar arguments see \cite{ip08,georgmarmousi}.  Our summary also
includes remarks on the relation to shape derivatives, as well as
computational details.

To simplify, the parameters $\kappa_{\rm e}, \kappa_{\rm i}, \beta$ are taken 
to be constant, real and positive. 
General conditions on the parameters $\kappa_{\rm e}$, $\kappa_{\rm i}$,
$\beta$ guaranteeing existence and uniqueness of a solution $u\in
H^{1}_{\rm loc}(\mathbb{R}^{2})$ for this  problem can be found  in
\cite{costabelstephan, kleinmanmartin, kressroach}.
When $\partial \Omega_{\rm i} \in C^2$, $u$ is in $H^2(\Omega_{\rm i}) \cup
H^2(\Omega_{\rm e})$, which ensures continuity away from the interface.

Let us recall that given two Banach spaces $X$, $Y$ and  a function 
${\cal F}: D({\cal F}) \subset X \longrightarrow Y$, its Fr\'echet derivative
${\cal F}': X \longrightarrow Y$ is a linear bounded operator satisfying
${\cal F}(x+ \xi) = {\cal F}(x) + {\cal F}'(x) \xi +  o (\xi)$
for $\xi \in X$ as $\| \xi \|_X \rightarrow 0$, for any $x \in X$. In terms
of the directional Gateaux derivative ${\cal F}'(x) \xi =  D_\xi {\cal F}(x)$
with  $ D_\xi{\cal F}(x)
= {\rm lim}_{\tau \rightarrow 0} {{\cal F}(x+ \tau \xi) - {\cal F}(x)
\over \tau }. $
In our context, Fr\'echet derivatives can be characterized as
solutions of adequate boundary value problems.

\subsection{Fr\'echet derivative with respect to the domain}
\label{sec:domain}
We consider  variable domains $\Omega_{\rm i} = \Omega^t$, whose boundaries 
$\Gamma^t$  are generated  from a smooth curve $ \Gamma^0 \in C^2$ (twice 
differentiable)  following a family of deformations  
$\Gamma^t = \left\{ \mathbf x + t \, \mathbf V(\mathbf x) \, | \, 
\mathbf x \in \Gamma^0 \right\}, $
along a smooth vector field $\mathbf V \in C^2 (\Gamma^0)$. The solutions of 
(\ref{forward}) with $\Omega_{\rm i}= \Omega^t$ are denoted by 
$u^t$. For small $t>0$, $ \Gamma^t \in C^2$ is a perturbation of $ \Gamma^0$. 
The deformation $\mathbf x^t = \phi^t(\mathbf x) = \mathbf x + t \, \mathbf 
V(\mathbf x)$ maps $\Omega^0$ to $\Omega^t$. For small $t$, $\phi^t$ is 
a diffeomorphism and its inverse $\eta^t$ maps $\Omega^t$ to $\Omega^0$.
We extend $\mathbf V$ to $\mathbb R^2$ in such a way that
it decays fast away from  $\Gamma^0$, while preserving the same regularity.

The operator ${\cal F}$ that assigns to $\mathbf V \in 
C^1(\partial \Omega_{\rm i})$ the far field values  $(u(\mathbf x_j))_{j=1}^N$ of the 
solution of the forward problem (\ref{forward}) is Fr\'echet differentiable
with derivative $\langle{\cal F}'(\partial \Omega_{\rm i}), \mathbf V \rangle = 
(v(\mathbf x_j))_{j=1}^N,$ where $v$ is the solution of  
\begin{equation} \hskip -6mm
\begin{array}{l}
\Delta v + \kappa_{\rm e}^2 v = 0   \quad
\mbox{\rm in $\Omega_{\rm e}$}, \qquad
\Delta v + \kappa_{\rm i}^2 v =0  \quad
\mbox{\rm in $\Omega_{\rm i}$},  \\ 
v^{-} \!-\! v^{+} \!=\!  - (\mathbf V \! \cdot \! \mathbf n)
\big( {\partial u^- \over \partial \mathbf n}
\!- \!{\partial u^+ \over \partial \mathbf n} \big)
\quad  \mbox{\rm on $\partial \Omega_{\rm i}$}, \\
\beta { \partial  v^{-}  \over \partial \mathbf n} 
\!-\! {\partial v^{+}  \over \partial \mathbf n} 
\!=\!  {d \over ds} \left[ (\mathbf V \! \cdot \! \mathbf n)
{d \over ds} (\beta u^- \!-\! u^+) \right] \!+\! (\mathbf V  \!\cdot\! \mathbf n)
(\beta \kappa_{\rm i}^2 u^- \!-\! \kappa_{\rm e}^2 u^+)
\quad \mbox{\rm on $\partial \Omega_{\rm i}$},  \\ 
{\rm lim}_{|\mathbf x| \rightarrow 0} |\mathbf x|^{1/2} 
\big({\partial \over \partial |\mathbf x|}v - \imath k_{\rm e} v \big) = 0,
\end{array} 
\label{frechet}
\end{equation}
and ${d \over ds}$ is the derivative with respect to the arclength,
as proven in \cite{hohage2d,hettlichfrechet}. 

In practice, the forward (\ref{forward}) and Fr\'echet 
(\ref{frechet}) problems are discretized and the solution operators they 
induce are represented by matrices. We use the boundary element
formulations introduced in \cite{javierbem1,javierbem2} to approximate
the solutions. To do so, we  recast (\ref{forward})  as a boundary
value problem for the scattered and the transmitted  fields with 
transmission conditions $u_{\rm tr}-u_{\rm sc}= u_{\rm inc}$ and $\beta
\partial_{\bf n} u_{\rm tr} - \partial_{\bf n} u_{\rm sc}=\partial_{\bf n} u_{\rm inc}$.

Note that the difference between two close domains 
$\Omega_{\nnu}$ and 
$\Omega_{\nnu + \boldsymbol \xi}$ parameterized by
(\ref{trigonometric1})-(\ref{trigonometric2})
defines a vector perturbation $\mathbf V$
of the boundary in terms of the difference
of their  parameterizations, which can be written in the form
with (\ref{trigonometric1})-(\ref{trigonometric2}) with
parameters $\boldsymbol \xi$ \cite{conca}.
Thus, the action of the Fr\'echet
derivative is identified with the action of a matrix on a vector, that is,
$\langle{\cal F}'(\partial \Omega_{\nnu}), \mathbf V \rangle $
becomes $\mathbf F(\nnu) \boldsymbol \xi 
\sim (v(\mathbf x_j))_{j=1}^N$. In this way, we compute the matrices
$ \mathbf F(\nnu)$ employed in Section \ref{sec:map}.

The characterization (\ref{frechet}) is proven using integral equations in 
\cite{hettlichfrechet,hohage2d} for $\mathbf V \in C^1$.
We give an alternative proof here using a variational approach inspired by 
the work  in \cite{dirichlet2D} for 2D exterior elasticity problems with zero 
Dirichlet boundary conditions on a moving boundary, which clarifies the 
role of the transmission conditions.

{\bf Theorem 1.} {\it Keeping the previous notations and assumptions, the 
`Fr\`echet derivative' of the far field of the solution $u$ of {\rm (\ref{forward})} 
with respect to the domain is given by the far field of the solution of the 
boundary value problem {\rm (\ref{frechet})}.} 

{\bf Proof.} The proof proceeds in the following  steps.

{\it Step 1: Variational formulation.} Firstly, we reformulate the 
transmission problem (\ref{forward})  as an equivalent boundary value 
problem posed in a bounded domain. Let $\Gamma_{R}$ be a large circle 
which encloses the objects $\Omega_{\rm i}$. 
The Dirichlet-to-Neumann  operator \cite{feijoooberai, kellergivoli}
associates to any Dirichlet data on $\Gamma_{R}$ the normal derivative 
of the solution of the exterior Dirichlet problem:
$L:H^{1/2}(\Gamma_{\rm R}) \longrightarrow  H^{-1/2}(\Gamma_{\rm R}),$
$f \longmapsto L(f) = \partial_{\bf n} v,$
where $v\in H^{1}_{\rm loc}(\mathbb{R}^{2} \setminus \overline{B}_{\rm R}), $
$B_{\rm R}:=B({\bf 0},R)$,  is the unique solution of
\begin{eqnarray*}
\left\{\begin{array}{l}  \Delta v+  \kappa_{\rm e}^2 v=0
\quad \mbox{in $\mathbb{R}^{2} \setminus \overline{B}_{\rm R}$}, \qquad
v=f \quad \mbox{on $\Gamma_{\rm R}$},\\ 
\displaystyle  \lim_{r\to\infty}r^{1/2}(\partial_r v-\imath
\kappa_{\rm e}
 v)=0.
\end{array}\right.
\end{eqnarray*}
In the sequel, $\partial_{\mathbf n}$ stands for the normal derivative at the 
interface and $\partial_{r}$ denotes the radial derivative. 
 $H^{1}_{\rm loc}(\mathbb{R}^{2} \setminus \overline{B}_{\rm R})$ is the
usual Sobolev space and $H^{1/2}(\Gamma_{\rm R})$ and $H^{-1/2}(\Gamma_{\rm R})$
denotes the standard trace spaces.  We replace (\ref{forward}) by an 
equivalent boundary value problem set in $B_{\rm R}$ with a non-reflecting 
boundary condition on $\Gamma_{\rm R}$ defined by the Dirichlet-to-Neuman 
operator:
\begin{equation}\label{transmisionBR}
\left\{\begin{array}{l}   \Delta u +  \kappa_{\rm e}^2 u=0 \quad
\mbox{in $\Omega'_{\rm e}:=B_{\rm R}\setminus \overline{\Omega}_{\rm i}$}, \qquad
\beta\Delta u +  \beta \kappa_{\rm i}^2 u=0 \quad \mbox
{in $\Omega_{\rm i}$},\\
u^--u^+=0 \quad \mbox{on $\partial \Omega_{\rm i}$}, \qquad
\beta\partial_{\bf n}u^--  \partial_{\bf n}u^+=0 \quad
\qquad\mbox{on $\partial \Omega_{\rm i}$},\\
\partial_{\bf n}(u-u_{\rm inc})=L(u-u_{\rm inc}) \quad
\mbox{on $\Gamma_{\rm R}$}.
\end{array} \right.
\end{equation}
Unit normals are exterior to $\partial \Omega_{\rm i}$ and to $\Gamma_{\rm R}$.
The solution $u$ of  (\ref{transmisionBR}) also solves  the
variational problem
\begin{eqnarray}\label{variational}
\begin{array}{l} u\in H^1(B_{\rm R}), \qquad
 b(\Omega_{\rm i}; u, w)= \ell(w),\qquad\forall w \in H^1(B_{\rm R}),
\end{array} \\
\begin{array}{rl}
 b(\Omega_{\rm i};u,w) &=  \ \int_{\Omega'_{\rm e}}
(  \nabla u\nabla {\overline w}- \kappa_{\rm e}^2 u{\overline w}) d\mathbf x
+\int_{\Omega_{\rm i}} (\beta\nabla u\nabla {\overline w}-
\beta \kappa_{\rm i}^2 u{\overline w}) d\mathbf x \\[1ex]
 & -\int_{\Gamma_{\rm R}}  Lu\,\overline{w} \, dS_{\mathbf x},
 \quad \forall u, w \in H^1(B_{\rm R}), \\[1ex]
\ell(w) & = \int_{\Gamma_{\rm R}} (\partial_{\bf n}u_{\rm inc}-Lu_{\rm inc})\,\overline{w}\, 
dS_{\mathbf x}, \quad \forall w \in H^1(B_{\rm R}).
\end{array} \nonumber
\end{eqnarray}

{\it Step 2: Deformed problems.}  Since $\mathbf V$ decreases
rapidly to zero away from $\Gamma^0$,  
$\phi^t(\Gamma_{\rm R})=\Gamma_{\rm R}$ and
$\phi^t(\mathbf x_j)=\mathbf x_j$, $j=1,\ldots,N$. Setting in
(\ref{variational})
$\Omega_{\rm i} =\phi^t(\Omega^0) = \Omega^t$ and
$\Omega_{\rm e}' =\phi^t(\Omega'_{\rm e}) = B_{\rm R}\setminus
\overline{\Omega}^t,$
the variational reformulations  in the deformed domains take the 
form: Find $u^t \in H^1(B_{\rm R})$ such that
\begin{eqnarray}\label{transformed} 
\begin{array}{ll} 
b^t(\Omega^t ; u^t, w)= \ell(w),\quad \forall w \in H^1(B_{\rm R}),
\end{array}   \\
\begin{array}{ll} 
b^t(\Omega^t;u,w)  &=  \int_{B_{\rm R}\setminus \overline{\Omega}^t }
(\nabla_{\mathbf x^t} u\nabla_{\mathbf x^t } {\overline w} -
\kappa_{\rm e}^2 u{\overline w}) d\mathbf x^t 
- \int_{\Gamma_{\rm R}} \!\! Lu\,\overline{w} dS_{\mathbf x} \\
&+  \int_{\Omega^t } (\beta\nabla_{\mathbf x^t }
u\nabla_{\mathbf x^t } {\overline w} - \beta \kappa_{\rm i}^2 u{\overline w}) 
d\mathbf x^t, \quad \forall u, w \in H^1(B_{\rm R}).
\end{array}  \nonumber
\end{eqnarray}

{\it Step 3: Change variables to initial configuration.} 
We now transform all the quantities appearing in (\ref{transformed}) back 
to the initial configuration $\Omega^0$. The process is similar
to transforming deformed configurations  back to a reference configuration 
in continuum mechanics.  
The deformation gradient is the jacobian of the change of variables 
\begin{eqnarray} \begin{array}{l}
\mathbf J^t(\mathbf x) = \nabla_{\mathbf x} \boldsymbol \phi^t(\mathbf x)  = 
\left({\partial x^t_i \over \partial x_j}(\mathbf x) \right) =
\mathbf I + t \, \nabla \mathbf V(\mathbf x), \end{array}
\label{jacobian}
\end{eqnarray} 
and its inverse $ (\mathbf J^t)^{-1}   
= \left({\partial x_i \over \partial x^t_j}\right)$ is the jacobian of the inverse
change of variables.
Then, volume and surface elements are related by
\begin{eqnarray} \begin{array}{l}
d \mathbf x^t = {\rm det} \, \mathbf J^t(\mathbf x) \, d \mathbf x, \quad
d S_{\mathbf x^t} = {\rm det} \, \mathbf J^t(\mathbf x)
\| \, (\mathbf J^t(\mathbf x))^{-T} \mathbf n \| dS_{\mathbf x},
\end{array} \label{changevolume}
\end{eqnarray}
and the chain rule for derivatives reads 
$\nabla_{\mathbf x} u(\mathbf x^t(\mathbf x)) = 
\mathbf (J^t(\mathbf x))^T \nabla_{\mathbf x^t} 
u(\mathbf x^t(\mathbf x)),$ that is,
$\nabla_{\mathbf x^t} u = (\mathbf J^t)^{-T} \nabla_{\mathbf x} u$.
For each component we have
\begin{eqnarray} \begin{array}{l}
{\partial u \over \partial x_\alpha^t}(\mathbf x^t(\mathbf x))  = 
{\partial u \over \partial x_k}(\mathbf x^t(\mathbf x)) 
(J^t)^{-1}_{k\beta}(\mathbf x).
\end{array} \label{changederivative}
\end{eqnarray}
Defining $\hat {u}(\mathbf x)= u^t  \circ \phi^t (\mathbf x)
= u^t (\mathbf x^t(\mathbf x))$, changing variables and using
(\ref{changevolume})-(\ref{changederivative}) yields:
\begin{eqnarray} \hskip -2.5cm \begin{array}{l}
b^t_{\rm i}(\Omega^t; u^t,w)   = \int_{\Omega^t}    \left[ \beta
{\partial u^t\over \partial x_{\alpha}^t }(\mathbf x^t) {\partial
\overline w \over \partial x_{\alpha}^t  }(\mathbf x^t)  - \beta \kappa_{\rm i}^2 
u^t(\mathbf x^t)  \overline w(\mathbf x^t) \right] \, d\mathbf x^t =  \\ 
\int_{\Omega^0} \beta \left[ {\partial \hat u \over \partial x_{p}}(\mathbf x) 
(J^t)^{-1}_{p \alpha}(\mathbf x)
{\partial \hat{\overline w} \over \partial x_{q} }(\mathbf x) 
(J^t)^{-1}_{q \alpha}(\mathbf x)  \, 
\!-\! \beta \kappa_{\rm i}^2   \hat u(\mathbf x)  \hat {\overline w}(\mathbf x)
\right] {\rm det} \, \mathbf J^t(\mathbf x) \, d \mathbf x  
= \hat b^t_{\rm i}(\Omega^0; \hat{u},\hat{w}).
\end{array} \label{changeb} 
\end{eqnarray}
A similar relation holds on $B_{\rm R}\setminus \overline{\Omega}^t$ defining 
$ b^t_{\rm e}(B_{\rm R}\setminus \overline{\Omega}^t; u^t,w) = 
\hat b^t_{\rm e}(B_{\rm R}\setminus \overline{\Omega}^0; \hat{u},\hat{w}).$
For $w \in H^1(B_{\rm R})$, we have $\hat {w} \in H^1(B_{\rm R})$. 
Therefore, we obtain the equivalent variational formulation:
Find $\hat{u} \in H^1(B_{\rm R})$ such that
\begin{equation*} \label{Variationalehat}
\hat b^t(\Omega^0; \hat{u}, w )  =
\displaystyle  \hat b^t_{\rm i}(\Omega^0; \hat{u}, w) +
\hat b^t_{\rm e}(B_{\rm R}\setminus \overline{\Omega}^0; \hat{u}, w) 
- \int_{\Gamma_{\rm R}}  L\hat u\,\overline{w} \, dS_{\mathbf x} = \ell(w),
\end{equation*}
for $w \in H^1(B_{\rm R})$.
Let us analyze the dependence on $t$ of the terms appearing in 
$\hat b^t_{\rm i}$ and $\hat b^t_{\rm e}$. From the definitions of the Jacobian matrices (\ref{jacobian}) we obtain \cite{feijoooberai,dirichlet2D}
\begin{eqnarray}   \begin{array}{l}
{\rm det} \, \mathbf J^t(\mathbf x)  \!=\! 1 \!+\! t \, {\rm div}(\mathbf V(\mathbf x) ) 
\!+\! O(t^2),  \quad
(\mathbf J^t)^{-1}(\mathbf x)  \!=\! \mathbf I \!-\! t \, \nabla \mathbf V(\mathbf x)  
\!+\! O(t^2). 
\end{array} \label{expandsurf}
\end{eqnarray}
Inserting (\ref{expandsurf}) in (\ref{changeb}) we find 
\begin{eqnarray} \begin{array}{l}
\ell(w) =\hat b^t(\Omega^0; \hat{u}, w)  \!=\!
b^0(\Omega^0; \hat{u}, w ) \!+\!
t[I_1(\hat{u}) \!+\! I_2(\hat{u}) \!+\! I_3(\hat{u})]
\!+\! O(t^2), \end{array}
\label{expandb}
\end{eqnarray}
where
\begin{eqnarray} \hskip -2cm
\begin{array}{ll}
 b^0(\Omega^0;\hat u,w) &= \int_{\Omega^0}   \left[ \beta
{\partial \hat u \over \partial x_{\alpha}  }
{\partial \overline w \over \partial x_{\alpha} }
- \beta \kappa_{\rm i}^2 \hat u {\overline w} \right] \, d \mathbf x +
 \int_{B_{\rm R}\setminus \overline{\Omega}^0}   \left[
{\partial \hat u \over \partial x_{\alpha}  }
{\partial \overline w \over \partial x_{\alpha} }
 - \kappa_{\rm e}^2 \hat u {\overline w} \right] \, d \mathbf x \\[1ex]
& -\int_{\Gamma_{\rm R}} L\hat u\,\overline{w} dS_{\mathbf x}, 
\end{array} \label{b0} \\ \hskip -2cm
\begin{array}{ll}
I_1(\hat{u}) =& \int_{\Omega^0}   \left[ \beta
{\partial \hat u \over \partial x_{\alpha}  }
{\partial \overline w \over \partial x_{\alpha} }
- \beta \kappa_{\rm i}^2 \hat u {\overline w} \right]
\, {\partial V_p \over \partial x_p}  \, d \mathbf x +
 \int_{B_{\rm R}\setminus \overline{\Omega}^0}   \left[
{\partial \hat u \over \partial x_{\alpha}  }
{\partial \overline w \over \partial x_{\alpha} }
- \kappa_{\rm e}^2 \hat u {\overline w} \right]
\, {\partial V_p \over \partial x_p}    \, d \mathbf x, \\
I_2(\hat{u})  =& - \int_{\Omega^0}  \beta
{\partial \hat u \over \partial x_{p}  } 
{\partial \overline w  \over \partial x_\alpha }
{\partial V_{p} \over \partial x_\alpha } d \mathbf x
- \int_{B_{\rm R}\setminus \overline{\Omega}^0}  
{\partial \hat u \over \partial x_{p}  } 
{\partial \overline w  \over \partial x_\alpha }
{\partial V_{p} \over \partial x_\alpha } 
d \mathbf x,  \\
I_3(\hat{u}) =&  
-  \int_{\Omega^0}  \beta {\partial \hat u \over \partial x_\alpha} 
{\partial \overline w \over \partial x_q } 
{\partial V_q \over \partial x_\alpha} \, d \mathbf x  
-  \int_{B_{\rm R}\setminus \overline{\Omega}^0}  
{\partial \hat u \over \partial x_\alpha} 
{\partial \overline w \over \partial x_q } 
{\partial V_q \over \partial x_\alpha} \, d \mathbf x.
 \end{array} \label{expandIs}
\end{eqnarray}

{\it Step 4. Variational problem for the domain derivative $u'$.}
Let us compare the transformed function $\hat{u} $ and the solution 
$u^0$ of $b^0(\Omega^0;u^0, w) = \ell(w)$. Thanks to (\ref{expandb}),
for any $w \in H^1(B_{\rm R})$ we have
\begin{eqnarray} \begin{array}{l}
b^0(\Omega^0;\hat{u}- u^0, w) =
- t [I_1(\hat{u})+I_2(\hat{u})+I_3(\hat{u})].
\end{array} \label{expandvariational}
\end{eqnarray}
Well posedness of the variational problems (\ref{transmisionBR}) with 
respect to changes in domains $\Omega^t$  implies uniform bounds on 
the solutions for $t \in [0,T]$:
$\|u^t\|_{H^1(B_{\rm R})} \leq C(T)$, $\| \hat{u} \| _{H^1(B_{\rm R})} \leq C(T).$
The  right hand side in (\ref{expandvariational})  tends to zero as 
$t \rightarrow 0$.
Well posedness of the variational problem again implies $\hat{u}
\rightarrow u^0$ in $H^1(B_{\rm R})$.

Dividing  (\ref{expandvariational}) by $t$, we find
$ b^0(\Omega^0;{\hat{u}- u^0 \over t}, w) =
- [I_1(\hat{u})+I_2(\hat{u})+I_3(\hat{u})]. $
Then, the limit $\dot {u} = {\rm lim}_{t \rightarrow 0} {\hat{u} - 
u^0 \over t}$ satisfies
\begin{eqnarray} 
\begin{array}{l}
b^0(\Omega^0;\dot {u}, w) = - [I_1(u^0)\!+\!I_2(u^0)\!+\!I_3(u^0)].
\end{array}
\label{equdot} 
\end{eqnarray}
The function $\dot {u}$ is the so called `material derivative',
that is, $\dot {u} = {\partial u \over \partial t} + \mathbf V \cdot
\nabla u^0 $. The domain derivative is then 
$u' = \dot {u} -  \mathbf V \cdot \nabla u^0$. Then,
\begin{eqnarray} \begin{array}{l}
b^0(\Omega^0; u', w) = b^0(\Omega^0;\dot {u}, w)
- b^0(\Omega^0;\mathbf V \cdot \nabla u^0 , w),
\end{array} \label{equprime}
\end{eqnarray}
where  $b^0(\Omega^0; \mathbf V \cdot \nabla u^0, w) 
$ is obtained from (\ref{b0}) replacing $\hat u$ by
$  {\partial u^0\over \partial x_p }  V_{p}.$

{\it Step 5. Differential equation for  ${u}'$.}
We evaluate the terms in the right hand side of (\ref{equdot})
to calculate the right hand side in (\ref{equprime}) using
(\ref{expandIs}).
Recalling the equations for $u^0$ and that $\mathbf V$ vanishes at $\Gamma_{\rm R}$, (\ref{equprime}) becomes after integrating by 
parts 
\begin{eqnarray} \hskip -2cm
\begin{array}{l}
b^0(\Omega^0; u', w) =
\int_{\Gamma^0}  \left(\beta {\partial (u^0)^-\over  \partial x_\alpha} -
{\partial (u^0)^+\over  \partial x_\alpha} \right) 
{\partial \overline w \over \partial x_q }  V_q n_\alpha \, d S_\mathbf x  \\[1ex]
-\int_{\Gamma^0}  \left(\beta {\partial (u^0)^-\over  \partial x_\alpha} -
{\partial (u^0)^+\over  \partial x_\alpha} \right) 
{\partial \overline w  \over \partial x_\alpha } V_p n_p d \mathbf x 
+ \int_{\Gamma^0} (\beta \kappa_{\rm i}^2 (u^0)^- - \kappa_{\rm e}^2 (u^0)^+)  
\overline w  
V_q n_q  \, d \mathbf x. 
\end{array} \label{equprimevar}
\end{eqnarray}
Integrating by parts in  $b^0(\Omega^0; u', w)$
for $w$ vanishing on $\Gamma^0$
this identity yields  the equations (\ref{frechet}) for $u'$:
$\Delta u' + \kappa^2_{\rm e} u'=0$ $   
\mbox{ in $B_{\rm R}\setminus \overline{\Omega}^0$}$,
$ \beta\Delta u' + \beta \kappa_{\rm i}^2 u'=0 $
$ \mbox{ in $\Omega^0$}$,
as well as the radiating boundary condition at infinity.

{\it Step 6: Transmission conditions for $u'$.}
When $\beta=1$, $u^0$ and its derivatives $\nabla u^0$
are continuous across  $\Gamma^0.$ As a result, $u'$
is continuous across $\Gamma^0$ and it belongs to
$H^1_{\rm loc}(\mathbb R^2)$.  
Integrating by parts in (\ref{equprimevar}) we get
\begin{eqnarray*}\begin{array}{ll}
\beta {\partial (u')^- \over \partial \mathbf n} -
{\partial (u')^+ \over \partial \mathbf n} =
(\beta \kappa_{\rm i}^2 (u^0)^- - \kappa_{\rm e}^2 (u^0)^+)    
(\mathbf V \cdot \mathbf n) & \mbox{on } \Gamma_0.
\end{array} \end{eqnarray*}
When $\beta \neq 1$, $\nabla u^0$ is not longer continuous across  
$\Gamma^0.$ Denoting by $\mathbf t$ and $\mathbf n$ the unit tangent
and normal vectors on $\Gamma^0$, we have
$\nabla (u^0)^{\pm} = (\nabla (u^0)^{\pm}  \cdot \mathbf t) \mathbf t +
(\nabla (u^0)^{\pm}  \cdot \mathbf n) \mathbf n. $
The relation $u'= \dot u - \mathbf V \cdot \nabla u^0$ gives the jump at the 
interface:
\begin{eqnarray*} \begin{array}{l}
(u')^- - (u')^+ =  - (\mathbf V \cdot \mathbf n) \left({\partial (u^0)^- \over
\partial \mathbf n} - {\partial (u^0)^+ \over \partial \mathbf n} \right).
\end{array} \label{equprimet1} \end{eqnarray*}
Notice that $u^0$ being continuous accross $\Gamma^0$, the tangent
derivatives ${\partial u^0\over \partial \mathbf t} =
\nabla u^0  \cdot \mathbf t$ too.
To obtain a transmission condition for the derivatives of
$u'$ at the interface  $\Gamma^0$ we revisit (\ref{equprimevar}).
The first term vanishes due to the transmission boundary conditions
satisfied by $u^0$ at $ \Gamma_0$. The second term can be rewritten
as
\begin{eqnarray*} \begin{array}{l}
-\int_{\Gamma^0}  \left(\beta {\partial (u^0)^-\over \partial \mathbf n} -
{\partial (u^0)^+\over \partial  \mathbf n} \right) 
{\partial \overline w  \over \partial \mathbf n} V_p n_p d \mathbf x
-\int_{\Gamma^0}  \left(\beta {\partial (u^0)^-\over  \partial \mathbf t} -
{\partial (u^0)^+\over  \partial \mathbf t} \right) 
{\partial \overline w  \over \partial \mathbf t} V_p n_p d \mathbf x.
\end{array} \end{eqnarray*}
The first integral vanished again due to the transmission boundary
conditions satisfied by $u^0$, whereas the second one, together with
the third integral on the right hand side of (\ref{equprimevar})
provides the transmission boundary condition in (\ref{frechet}).
\proofend

Given the characterization of the Fr\'echet derivative (\ref{frechet}), we
find an expression for the shape derivative of the cost (\ref{cost})  along a 
smooth vector field ${\mathbf V}$, defined as $DJ_c(\Omega_{\rm i})\cdot{\bf V}=
\frac{d}{dt}\,J_c(\phi^t(\Omega_{\rm i})) \big|_{t=0},$ keeping the previous
notations and following \cite{sims18,jcp19}:
\begin{eqnarray}\label{shapederivative}    \begin{array}{ll}
\langle DJ_c({\mathbb{R}}^2\setminus \overline{ \Omega_{\rm i}}), \mathbf 
V \rangle   = &  
\mbox{\rm Re} \big[ \int_{\partial \Omega_{\rm i}} \! \big( (1 - \beta) \big(\beta  
\partial_{\bf n} u^- \partial_{\bf n} {\overline{p}}^- + \partial_{\bf t} u^- 
\partial_{\bf t} {\overline{p}}^- \big)   \\[1ex]
 & + (\beta \kappa_{\rm i}^2 - \kappa_{\rm e}^2  )\, u \,{\overline{p}}    
\big) \mathbf V \cdot \mathbf n \,dS_{\mathbf x} \big],
\end{array} \end{eqnarray}
where $u$ and $\overline{p}$ solve the forward and adjoint problems 
(\ref{forward}) and  (\ref{adjoint})
\begin{eqnarray}
\left|\begin{array}{l} \Delta \overline{p} + \kappa_{\rm e}^2 \overline{p}= 
\sum_{j=1}^N \chi(\mathbf x_j)\delta_{\mathbf x_j} \quad
\mbox{\rm in ${\mathbb{R}}^2 \setminus \overline{ \Omega}_{i}$}, \quad
\Delta \overline{p} + \kappa_{\rm i}^2 \overline{p}=0
\quad \mbox{\rm in $\Omega_{\rm i}$},\\ 
\overline{p}^--\overline{p}^+=0 \quad \mbox{\rm on $\partial \Omega_{\rm i}$},  
\quad \beta \partial_{\bf n}\overline{p}^-- 
\partial_{\bf n}\overline{p}^+=0 \quad \mbox{\rm on
$\partial \Omega_{\rm i}$},\\ 
\displaystyle \lim_{r\to\infty}r^{1/2}\left(\partial_r \overline{p}- \imath
k_{\rm e} \overline{p} \right)=0,\end{array} \right. 
\label{adjoint} \\
\chi(u(\mathbf x_j)) = \left\{ \begin{array}{l}
\overline{\dnj-u(\mathbf x_j)}, \quad 
\mbox{\rm for  } f(u)=u, \\
2 (\dnj -|u(\mathbf x_j)|^2) \overline{u(\mathbf x_j)}, 
\quad \mbox{\rm for  } f(u)=|u|^2,
\end{array} \right. \label{def_d}
\end{eqnarray} 
with $r = |\mathbf x|.$
This expression agrees with that established in \cite{ln08,sims18} when $\mathbf V = V_{\rm n} \mathbf n$. The same proofs hold for general fields $\mathbf V$ keeping track of the terms involving tangential components and the transmission boundary conditions.
Writing $\overline p = \overline p' + \overline{p}_{\rm inc}$, where $\overline{p}_{\rm inc}$  is given by (\ref{pinc}),  $\overline p' $ is a solution of a problem of the form (\ref{forward}) with incident wave $\overline{p}_{\rm inc}$, 
solvable by boundary elements too. 

For the functional (\ref{regcost}), the derivatives with respect to $\nnu$  are given in Section \ref{sec:map} in terms of Fr\'echet derivatives. In terms of
shape derivatives, the gradient method reads
$ \nu_{n+1,m}^\ell  = \nu_{n,m}^\ell - \tau_{n,m}^\ell  \big[ \sigma_{\rm noise}^{-2}
\langle DJ_c(\mathbb R^2 \setminus \overline{\Omega}_{\nnu_n}),
\mathbf V_{m} \rangle $ $+ \sigma_m^{-2}
(\nnu_{n,m} - \nnu_{0,m}) \big],  $
for $\ell=1,\ldots,L$, $m=1,\ldots,2M+3$, $\tau_{n,m}^\ell>0$ small, with directions
$\mathbf V_1=(1,0),$  $\mathbf V_2=(0,1),$  
$\mathbf V_{2+m+1}= \cos(2 \pi m t) (\cos(2 \pi t), \sin(2 \pi t) ), $
for $ m=0,\ldots,M $ and
$\mathbf V_{3+M+m}= \sin(2 \pi m t) (\cos(2 \pi t), \sin(2 \pi t) ), $ for
$m=1,\ldots,M.$
We have used the above formulas to check codes and iterative procedures.

\subsection{Fr\'echet derivative with respect to the coefficients}
\label{sec:coefficients}

An analogous but much simpler procedure to that followed the previous
section, provides an expression for Fr\'echet derivatives when
$\kappa_{\rm i}$ is a constant.

{\bf Theorem 2.} {\it Keeping the notations and assumptions of the
previous section, the `Fr\`echet derivative' of the far field of the solution $u$ 
of {\rm (\ref{forward})} with respect to $\kappa_{\rm i}$ is given by the far field of 
the solution of the boundary value problem {\rm (\ref{frechetc})}:
\begin{equation} 
\begin{array}{l}
\Delta v + \kappa_{\rm e}^2 v = 0   \quad
\mbox{\rm in $\Omega_{\rm e}$}, \qquad
\Delta v + \kappa_{\rm i}^2 v = - 2 \kappa_{\rm i} u \quad 
\mbox{\rm in $\Omega_{\rm i}$},  \\ 
v^{-} \!-\! v^{+} \!=\!  0
\quad  \mbox{\rm on $\partial \Omega_{\rm i}$},  \quad
\beta { \partial  v^{-}  \over \partial \mathbf n} 
\!-\! {\partial v^{+}  \over \partial \mathbf n} 
\!=\!  0 \quad \mbox{\rm on $\partial \Omega_{\rm i}$},  \\ 
{\rm lim}_{|\mathbf x| \rightarrow 0} |\mathbf x|^{1/2} 
\big({\partial \over \partial |\mathbf x|}v - \imath k_{\rm e} v \big) = 0.
\end{array} 
\label{frechetc}
\end{equation}}
{\bf Proof.} Let $u$ be the solution of (\ref{forward}), $\kappa_{\rm i}$
constant, and $u^t$ the solution of (\ref{forward}) with coefficient
$\kappa_{\rm i}^t= \kappa_{\rm i} + t$. Then $u^t-u$ is a solution of
(\ref{forward}) with right hand side $f^t= - 2 \kappa_{\rm i} t u^t - t^2 u^t$
and zero incident wave.
As argued in Step 1 of the proof of Theorem 1, the forward
problem (\ref{forward}) admits the variational reformulation
 (\ref{transmisionBR}).
Well-posedness of the variational problems (\ref{transmisionBR}) with 
respect to changes in the coefficients  implies uniform bounds on 
the solutions for $t \in [0,T]$: $\|u^t\|_{H^1(B_{\rm R})} \leq C(T)$.
Thus, the  right hand sides $f^t$ tend to zero as $t \rightarrow 0$.
Well-posedness of the variational problem again implies that $u^t
\rightarrow u^0=u$ in $H^1(B_{\rm R})$. 
The quotients ${u^t - u^0 \over t}$ are solutions of (\ref{forward})
with right hand side ${f^t - f^0 \over t} = - 2 \kappa_{\rm i}  u^t - t u^t$,
which tends to $ - 2 \kappa_{\rm i}  u $ in $H^1(B_{\rm R})$. 
Well-posedness of the variational problem again with respect to
the right hand side implies that the limit 
$\dot {u} = {\rm lim}_{t \rightarrow 0} {u^t - u^0 \over t}$ is a solution
of (\ref{frechetc}). \proofend

This boundary problem has a non zero right hand side. We may solve it by 
finite elements. However, in an iterative optimization procedure that requires 
solving problems of the form (\ref{frechetc}) for different $\Omega_{\rm i}$
in each iteration, we have to use fine meshes and remesh each time 
we change the proposed objects, which is expensive.
Alternatively, we can solve for the right hand side in the whole space by convolution with the Green function of the exterior Helmholtz problem and correct it solving a transmission problem with zero source by BEM. The convolution makes this option equally expensive. From the computational 
point of view it is more convenient to use the definition 
${\partial u(\kappa_{\rm i}) \over \partial t} = {u(\kappa_{\rm i}+t) + u(\kappa_{\rm i}) \over t}$
for small $t>0$ to approximate it. In this way we only have to solve an
additional forward problem with parameter $\kappa_{\rm i}+t$ by BEM per 
iteration. In this way we complete the Fr\'echet matrices $\mathbf 
F(\nnu)$ with an additional column to obtain the Fr\'echet
matrices $\mathbf F(\nnu, \kappa_{\rm i}).$

\end{document}